\documentclass{amsart}
\usepackage{amsmath,amssymb,amsxtra,amscd}
\usepackage[all]{xy}
\usepackage{verbatim}
\usepackage[vcentermath,enableskew]{rcyoungtab}
\usepackage{pstricks}
\usepackage{pst-node}

\begin{document}

\newcommand{\ul}[1]{\underline{#1}}

\Yboxdim{13pt}\Yinterspace{1pt}

\newcommand{\bl}{{}}

\newcommand{\ba}{\bar{1}}
\newcommand{\bb}{\bar{2}}
\newcommand{\bc}{\bar{3}}
\newcommand{\bd}{\bar{4}}

\newcommand{\ca}{\check{1}}
\newcommand{\cb}{\check{2}}
\newcommand{\cc}{\check{3}}
\newcommand{\cd}{\check{4}}
\newcommand{\ce}{\check{5}}
\newcommand{\cf}{\check{6}}
\newcommand{\cg}{\check{7}}
\newcommand{\ch}{\check{8}}
\newcommand{\ci}{\check{9}}
\newcommand{\cj}{\check{a}}
\newcommand{\ck}{\check{b}}
\newcommand{\cl}{\check{c}}
\newcommand{\cm}{\check{d}}

\newcommand{\ha}{\hat{1}}
\newcommand{\hb}{\hat{2}}
\newcommand{\hc}{\hat{3}}
\newcommand{\hd}{\hat{4}}
\newcommand{\he}{\hat{5}}
\newcommand{\hf}{\hat{6}}
\newcommand{\hg}{\hat{7}}
\newcommand{\hh}{\hat{8}}
\newcommand{\hi}{\hat{9}}
\newcommand{\hj}{\hat{a}}
\newcommand{\hk}{\hat{b}}
\newcommand{\hl}{\hat{c}}
\newcommand{\hm}{\hat{d}}

\newcommand{\sta}{1^*}
\newcommand{\stb}{2^*}
\newcommand{\stc}{3^*}
\newcommand{\std}{4^*}
\newcommand{\ste}{5^*}
\newcommand{\stf}{6^*}
\newcommand{\stg}{7^*}
\newcommand{\sth}{8^*}
\newcommand{\sti}{9^*}
\newcommand{\stj}{a^*}
\newcommand{\stk}{b^*}
\newcommand{\stl}{c^*}
\newcommand{\stm}{d^*}

\newcommand{\pmt}[2]{\begin{pmatrix} #1 \\ #2 \end{pmatrix}}%
\newcommand{\twocol}[1]{\multicolumn{2}{c}{#1}}

\newcommand{\shape}{\mathrm{shape}}
\newcommand{\dt}{\cdot}

\numberwithin{equation}{section}

\newtheorem{theorem}{Theorem}
\newtheorem{prop}[theorem]{Proposition}
\newtheorem{conjecture}[theorem]{Conjecture}
\newtheorem{lemma}[theorem]{Lemma}
\newtheorem{cor}[theorem]{Corollary}
\theoremstyle{definition}
\newtheorem{remark}[theorem]{Remark}
\newtheorem{question}[theorem]{Question}
\newtheorem{example}[theorem]{Example}
\newtheorem{assumption}[theorem]{Assumption}

\newcommand{\anti}{\#}
\newcommand{\bilet}[2]{\begin{pmatrix} #2 \\ #1 \end{pmatrix}}
\newcommand{\CA}{\mathcal{C}'}
\newcommand{\CC}{\mathcal{C}}
\newcommand{\cck}{\check{c}}
\newcommand{\cht}{\widehat{c}}
\newcommand{\cocharge}{\overline{c}}
\newcommand{\coenergy}{\Db}
\newcommand{\Db}{\overline{D}}
\newcommand{\dck}{\check{d}}
\newcommand{\dckh}{\widehat{\dck}}
\newcommand{\dcm}{\check{d}_-}
\newcommand{\dcp}{\check{d}_+}
\newcommand{\eh}{\hat{e}}
\newcommand{\ev}{\mathrm{e}}
\newcommand{\fh}{\hat{f}}
\newcommand{\gb}{\overline{\gf}}
\newcommand{\gf}{\mathfrak{g}}
\newcommand{\Hb}{\overline{H}}
\newcommand{\Image}{\mathrm{Image}}
\newcommand{\Kb}{\overline{K}}
\newcommand{\KK}{f}
\newcommand{\Knuth}{\sim}
\newcommand{\La}{\Lambda}
\newcommand{\la}{\lambda}
\newcommand{\Lab}{\overline{\La}}
\newcommand{\lad}{\la^{\bullet}}
\newcommand{\LR}{\mathrm{LR}}
\newcommand{\Mb}{\overline{M}}
\newcommand{\nn}{\varnothing}
\newcommand{\nuh}{{\hat{\nu}}}
\newcommand{\Par}{\mathcal{P}}
\newcommand{\Pb}{\overline{P}}
\newcommand{\Pcl}{P_{cl}}
\newcommand{\Pns}{P^{\mathrm{ns}}}
\newcommand{\Pt}{\tilde{P}}
\newcommand{\Pv}{P^v}
\newcommand{\Rv}{R^v}
\newcommand{\Qt}{\tilde{Q}}
\newcommand{\rs}{\mathrm{rs}}
\newcommand{\rsh}{\hat{\rs}}
\newcommand{\se}{\diamondsuit}
\newcommand{\swch}{s}
\newcommand{\thet}[3]{{}^{#1}\theta_{#2}^{#3}}
\newcommand{\Tt}{\tilde{T}}
\newcommand{\Vh}{\hat{V}}
\newcommand{\ve}{\varepsilon}
\newcommand{\vn}{\varnothing}
\newcommand{\vp}{\varphi}
\newcommand{\VX}{VX}
\newcommand{\word}{\mathrm{word}}
\newcommand{\wt}{\mathrm{wt}}
\newcommand{\Xb}{\overline{X}}
\newcommand{\Z}{\mathbb{Z}}
\newcommand{\zv}{z^\vee}
\newcommand{\inner}[2]{\langle#1,#2\rangle}

\newcommand{\hdom}{{\begin{picture}(8,4)
\multiput(0,0)(4,0){3}{\line(0,1){4}}%
\multiput(0,0)(0,4){2}{\line(1,0){8}}\end{picture}}}

\newcommand{\vdom}{{\begin{picture}(4,8)
\multiput(0,0)(0,4){3}{\line(1,0){4}}%
\multiput(0,0)(4,0){2}{\line(0,1){8}}\end{picture}}}

\newcommand{\cell}{{\begin{picture}(4,4)
\multiput(0,0)(4,0){2}{\line(0,1){4}}%
\multiput(0,0)(0,4){2}{\line(1,0){4}}\end{picture}}}

\newcommand{\theset}{\{\nn,\,\cell\,,\,\hdom\,,\,\vdom\,\}}
\newcommand{\theneset}{\{\cell\,,\,\hdom\,,\,\vdom\,\}}
\newcommand{\ourset}{\{\nn,\cell\,,\hdom\,\}}
\newcommand{\onetwoset}{\{\cell\,,\hdom\,\}}

\title{On the $X=M=K$ conjecture}
\author{Mark Shimozono}
\address{Department of Mathematics \\
Virginia Tech \\
Blacksburg, VA 24061-0123 USA}
\thanks{Supported in part by NSF grant DMS-0401012}

\begin{abstract} In the large rank limit,
for any nonexceptional affine algebra, the graded branching
multiplicities known as one-dimensional sums, are conjectured to
have a simple relationship with those of type $A$, which are known
as generalized Kostka polynomials. This is called the $X=M=K$
conjecture. It is proved for tensor products of the ``symmetric
power" Kirillov-Reshetikhin modules for all nonexceptional affine
algebras except those whose Dynkin diagrams are isomorphic to that
of untwisted affine type $D$ near the zero node. Combined with
results of Lecouvey, this realizes the above one-dimensional sums of
affine type $C$, as affine Kazhdan-Lusztig polynomials (and
conjecturally for type $D$).
\end{abstract}

\maketitle

\section{Introduction}

\subsection{$X=M$ conjecture} \label{ss:X=M} Motivated by the
study of two-dimensional solvable lattice models, the seminal paper
\cite{KMN92} introduced the path model for the crystal graph of an
irreducible integrable highest weight module over a quantum affine
algebra $U'_q(\gf)$. A path is a semi-infinite sequence of elements
taken from the crystal graph of a suitable finite-dimensional
$U'_q(\gf)$-module. To apply the path model these suitable modules
and crystal graphs must be constructed. Unlike irreducible
integrable highest weight modules over a quantum Kac-Moody algebra,
it is rare for a finite-dimensional $U'_q(\gf)$-module to have a
crystal base. A number of suitable modules and crystal graphs had
been constructed in the literature
\cite{JMO00,KKM94,KMN92a,Ko99,Y98}.

Inspired by the work of Kirillov and Reshetikhin \cite{KR90} on
finite-dimensional representations of Yangians, the papers
\cite{HKOTT01,HKOTY99} conjectured the existence of a suitable
family of finite-dimensional $U'_q(\gf)$-modules $\{W_s^{(r)} \}$
called Kirillov-Reshetikhin (KR) modules. The module $W_s^{(r)}$ is
conjectured to be irreducible and to have a crystal base $B^{r,s}$.
The KR modules are indexed by pairs $(r,s)$ where $r$ is a node of
the Dynkin diagram of $\gb$ and $s\in \Z_{>0}$. The family of KR
modules is complete in the following sense: it is expected that
every finite-dimensional irreducible $U'_q(\gf)$-module with affine
crystal base, is a tensor product of KR modules \cite{Kas03}.

Let $W^L = \bigotimes (W_s^{(r)})^{\otimes L_s^{(r)}}$ be a finite
tensor product of KR modules, where $L=(L_s^{(r)})$ is a collection
of nonnegative integers. The tensor product $W^L$ has a
$U_q(\gb)$-equivariant grading by the coenergy function
\cite{HKOTT01,KMN92,NY97}. The \textbf{one-dimensional sum}
$\Xb_{L,\la}(t)$ is by definition the graded multiplicity in $W^L$
of the irreducible $U_q(\gb)$-module $V(\la)$ of highest weight
$\la$. The one-dimensional sum may be defined solely in terms of the
combinatorics of the KR crystals $B^{r,s}$. For affine type A, the
one-dimensional sums are known as \textbf{generalized Kostka
polynomials} \cite{NY97,ScWa99,S02} and are well-understood
combinatorially.

Based on considerations arising from the Bethe Ansatz \cite{KR90},
it was conjectured in \cite{HKOTT01,HKOTY99} that the
one-dimensional sum $\Xb_{L,\la}(t)$ is equal to a specific sum of
products of $t$-binomial coefficients, known as the
\textbf{fermionic formula} $\Mb_{L,\la}(t)$. This is the
\textbf{$\mathbf{X=M}$ conjecture}. The fermionic formula may be
expressed as a weighted sum over a set of combinatorial objects
called rigged configurations. For untwisted affine type $A$, $X=M$
was proved by exhibiting a grade-preserving bijection from classical
highest weight vectors in the affine crystal graph, to rigged
configurations \cite{KKR88,KR88,KSS02}. The $X=M$ conjecture has
also been proved for tensor products of the ``symmetric power" KR
crystal $B^{1,s}$ for nonexceptional affine algebras \cite{ScS04}
based partially on previous work for tensor powers of the ``vector"
KR crystal $B^{1,1}$ \cite{OSS03a}. $X=M$ has also been proved for
tensor products of the ``exterior power" KR crystals (those of the
form $B^{r,1}$) for the root systems $C_n^{(1)}$, $A_{2n-2}^{(2)}$,
$D_{n+1}^{(2)}$ \cite{OSS03} and $D_n^{(1)}$ \cite{Sci04}. In the
ungraded case $t=1$, $X=M$ is known in many cases (for signed
binomial coefficients in the $M$ formula) by combining \cite{Ch01}
and \cite{HKOTT01,HKOTY99}.

\subsection{Large rank}
\label{ss:limit} We show that for large rank the $X$ or $M$
polynomials have a surprisingly simple expression (the $K$ formula)
that involves affine type $A$ only.

Let $\{\gf_n\}$ be a family of nonexceptional affine algebras. They
are depicted in Figure \ref{fig:rootsys}. Let $\gb_n$ be the simple
Lie subalgebra of $\gf_n$ whose Dynkin diagram is obtained by
removing the zero node. The subscript $n$ indicates that $\gb_n$ has
rank $n$. Fix a partition $\la=(\la_1,\la_2,\dotsc,\la_r)$. For $n$
large, the partition $\la$ may be identified with the dominant
$\gb_n$-weight $\sum_i (\la_i-\la_{i+1}) \omega_i$, where $\omega_i$
is the $i$-th fundamental weight of $\gb_n$. Since the rank is
assumed to be large with respect to $\la$, $\la$ involves no spin
weights.

\begin{figure}
\begin{tabular}[t]{|c|c|c|} \hline
$\se$ & $X_N^{(r)}$ & \text{Dynkin diagram} \\ \hline $\nn$ &
\begin{minipage}{4em}
\begin{center}
$A_n^{(1)}$ \\ $(n \ge 2)$
\end{center}
\end{minipage} &
\raisebox{-7mm}{%
\scalebox{.9}{%
\psset{unit=1cm}%
\pspicture(-.25,-.5)(5.25,1.35)%
\rput(0,0){\circlenode{A1}{}}%
\rput(1,0){\circlenode{A2}{}}%
\rput(1.75,0){\pnode{A3}}%
\rput(2.5,0){$\dotsm\dotsm$}%
\rput(3.25,0){\pnode{A4}}%
\rput(4,0){\circlenode{A5}{}}%
\rput(5,0){\circlenode{A6}{}}%
\rput(2.5,.75){\circlenode{B}{}}%
\ncline{A1}{A2}%
\ncline{A2}{A3}%
\ncline{A4}{A5}%
\ncline{A5}{A6}%
\ncline{A1}{B}%
\ncline{A6}{B}%
\rput(0,-.3){1}%
\rput(1,-.3){2}%
\rput(4,-.3){n-1}%
\rput(5,-.3){n}%
\rput(2.5,1.1){0}%
\endpspicture}}%
\\ \hline %
$\cell$ &
\begin{minipage}{4em}
\begin{center}
$D_{n+1}^{(2)}$\\$(n \ge 2)$
\end{center}
\end{minipage}&
\raisebox{-4.5mm}{%
\scalebox{.9}{%
\psset{unit=1cm}%
\pspicture(-.25,-.6)(6.25,.6)%
\rput(0,0){\circlenode{A0}{}}%
\rput(1,0){\circlenode{A1}{}}%
\rput(2,0){\circlenode{A2}{}}%
\rput(2.75,0){\pnode{A3}}%
\rput(3.5,0){$\dotsm\dotsm$}%
\rput(4.25,0){\pnode{A4}}%
\rput(5,0){\circlenode{A5}{}}%
\rput(6,0){\circlenode{A6}{}}%
\ncline[doubleline=true]{<-}{A0}{A1}%
\ncline{A1}{A2}%
\ncline{A2}{A3}%
\ncline{A4}{A5}%
\ncline[doubleline=true]{->}{A5}{A6}%
\rput(0,-.3){0}%
\rput(1,-.3){1}%
\rput(2,-.3){2}%
\rput(5,-.3){n-1}%
\rput(6,-.3){n}%
\endpspicture}}
\\ \hline
$\cell$ &
\begin{minipage}{4em}
\begin{center}
$A_{2n}^{(2)}$ \\ $(n \ge 2)$
\end{center}
\end{minipage}&
\raisebox{-4.5mm}{%
\scalebox{.9}{%
\psset{unit=1cm}%
\pspicture(-.25,-.6)(6.25,.6)%
\rput(0,0){\circlenode{A0}{}}%
\rput(1,0){\circlenode{A1}{}}%
\rput(2,0){\circlenode{A2}{}}%
\rput(2.75,0){\pnode{A3}}%
\rput(3.5,0){$\dotsm\dotsm$}%
\rput(4.25,0){\pnode{A4}}%
\rput(5,0){\circlenode{A5}{}}%
\rput(6,0){\circlenode{A6}{}}%
\ncline[doubleline=true]{<-}{A0}{A1}%
\ncline{A1}{A2}%
\ncline{A2}{A3}%
\ncline{A4}{A5}%
\ncline[doubleline=true]{<-}{A5}{A6}%
\rput(0,-.3){0}%
\rput(1,-.3){1}%
\rput(2,-.3){2}%
\rput(5,-.3){n-1}%
\rput(6,-.3){n}%
\endpspicture}}
\\ \hline
$\hdom$ &
\begin{minipage}{4em}
\begin{center}
$C_n^{(1)}$\\$(n \ge 2)$
\end{center}
\end{minipage}&
\raisebox{-4.5mm}{%
\scalebox{.9}{%
\psset{unit=1cm}%
\pspicture(-.25,-.6)(6.25,.6)%
\rput(0,0){\circlenode{A0}{}}%
\rput(1,0){\circlenode{A1}{}}%
\rput(2,0){\circlenode{A2}{}}%
\rput(2.75,0){\pnode{A3}}%
\rput(3.5,0){$\dotsm\dotsm$}%
\rput(4.25,0){\pnode{A4}}%
\rput(5,0){\circlenode{A5}{}}%
\rput(6,0){\circlenode{A6}{}}%
\ncline[doubleline=true]{->}{A0}{A1}%
\ncline{A1}{A2}%
\ncline{A2}{A3}%
\ncline{A4}{A5}%
\ncline[doubleline=true]{<-}{A5}{A6}%
\rput(0,-.3){0}%
\rput(1,-.3){1}%
\rput(2,-.3){2}%
\rput(5,-.3){n-1}%
\rput(6,-.3){n}%
\endpspicture}}
\\ \hline
$\hdom$ &
\begin{minipage}{4em}
\begin{center}
$A_{2n}^{(2)\dagger}$\\$(n \ge 2)$
\end{center}
\end{minipage}&
\raisebox{-4.5mm}{%
\scalebox{.9}{%
\psset{unit=1cm}%
\pspicture(-.25,-.6)(6.25,.6)%
\rput(0,0){\circlenode{A0}{}}%
\rput(1,0){\circlenode{A1}{}}%
\rput(2,0){\circlenode{A2}{}}%
\rput(2.75,0){\pnode{A3}}%
\rput(3.5,0){$\dotsm\dotsm$}%
\rput(4.25,0){\pnode{A4}}%
\rput(5,0){\circlenode{A5}{}}%
\rput(6,0){\circlenode{A6}{}}%
\ncline[doubleline=true]{->}{A0}{A1}%
\ncline{A1}{A2}%
\ncline{A2}{A3}%
\ncline{A4}{A5}%
\ncline[doubleline=true]{<-}{A5}{A6}%
\rput(0,-.3){0}%
\rput(1,-.3){1}%
\rput(2,-.3){2}%
\rput(5,-.3){n-1}%
\rput(6,-.3){n}%
\endpspicture}}
\\ \hline
$\vdom$ &
\begin{minipage}{4em}
\begin{center}
$B_n^{(1)}$\\$(n \ge 3)$
\end{center}
\end{minipage}&
\raisebox{-7mm}{%
\scalebox{.9}{%
\psset{unit=1cm}%
\pspicture(-.25,-.6)(6.25,1.3)%
\rput(0,0){\circlenode{A1}{}}%
\rput(1,0){\circlenode{A2}{}}%
\rput(1,1){\circlenode{A0}{}}%
\rput(2,0){\circlenode{A3}{}}%
\rput(2.75,0){\pnode{A4}}%
\rput(3.5,0){$\dotsm\dotsm$}%
\rput(4.25,0){\pnode{A5}}%
\rput(5,0){\circlenode{A6}{}}%
\rput(6,0){\circlenode{A7}{}}%
\ncline{A1}{A2}%
\ncline{A2}{A0}%
\ncline{A2}{A3}%
\ncline{A3}{A4}%
\ncline{A5}{A6}%
\ncline[doubleline=true]{->}{A6}{A7}%
\rput(0,-.3){1}%
\rput(1,-.3){2}%
\rput(1.3,1){0}%
\rput(2,-.3){3}%
\rput(5,-.3){n-1}%
\rput(6,-.3){n}%
\endpspicture}}
\\ \hline
$\vdom$ &
\begin{minipage}{4em}
\begin{center}
$D_n^{(1)}$\\$(n \ge 4)$
\end{center}
\end{minipage}&
\raisebox{-7mm}{%
\scalebox{.9}{%
\psset{unit=1cm}%
\pspicture(-.25,-.6)(6.25,1.3)%
\rput(0,0){\circlenode{A1}{}}%
\rput(1,0){\circlenode{A2}{}}%
\rput(1,1){\circlenode{A0}{}}%
\rput(2,0){\circlenode{A3}{}}%
\rput(2.75,0){\pnode{A4}}%
\rput(3.5,0){$\dotsm\dotsm$}%
\rput(4.25,0){\pnode{A5}}%
\rput(5,0){\circlenode{A6}{}}%
\rput(6,0){\circlenode{A7}{}}%
\rput(5,1){\circlenode{A8}{}}%
\ncline{A1}{A2}%
\ncline{A2}{A0}%
\ncline{A2}{A3}%
\ncline{A3}{A4}%
\ncline{A5}{A6}%
\ncline{A6}{A8}%
\ncline{A6}{A7}%
\rput(0,-.3){1}%
\rput(1,-.3){2}%
\rput(1.3,1){0}%
\rput(2,-.3){3}%
\rput(5,-.3){n-2}%
\rput(6,-.3){n-1}%
\rput(5.3,1){n}%
\endpspicture}}
\\ \hline
$\vdom$ &
\begin{minipage}{4em}
\begin{center}
$A_{2n-1}^{(2)}$\\$(n \ge 3)$
\end{center}
\end{minipage}&
\raisebox{-7mm}{%
\scalebox{.9}{%
\psset{unit=1cm}%
\pspicture(-.25,-.6)(6.25,1.3)%
\rput(0,0){\circlenode{A1}{}}%
\rput(1,0){\circlenode{A2}{}}%
\rput(1,1){\circlenode{A0}{}}%
\rput(2,0){\circlenode{A3}{}}%
\rput(2.75,0){\pnode{A4}}%
\rput(3.5,0){$\dotsm\dotsm$}%
\rput(4.25,0){\pnode{A5}}%
\rput(5,0){\circlenode{A6}{}}%
\rput(6,0){\circlenode{A7}{}}%
\ncline{A1}{A2}%
\ncline{A2}{A0}%
\ncline{A2}{A3}%
\ncline{A3}{A4}%
\ncline{A5}{A6}%
\ncline[doubleline=true]{<-}{A6}{A7}%
\rput(0,-.3){1}%
\rput(1,-.3){2}%
\rput(1.3,1){0}%
\rput(2,-.3){3}%
\rput(5,-.3){n-1}%
\rput(6,-.3){n}%
\endpspicture}}
\\ \hline
\end{tabular}
\caption{\label{fig:rootsys} Dynkin diagrams for nonexceptional
affine root systems}
\end{figure}

Let $L=(L^{(r)}_s)$ be a collection of nonnegative integers, only
finitely many of which are nonzero, indexed by $(r,s)\in\Z_{>0}^2$.
For $n$ large, one may associate to $L$ the tensor product of KR
modules $W^L=\bigotimes_{r,s} W_s^{(r) \otimes L^{(r)}_s}$.

\begin{lemma} \label{lem:afflimit} \cite{SZ04} Let $\la$ and $L=(L^{(r)}_s)$
be as above, and $\mathcal{F}=\{\gf_n\}$ a nonexceptional family
of affine algebras. Then there is a limiting polynomial
\begin{equation}
\Mb^{\mathcal{F}}_{L,\la}(t) =
  \lim_{n\rightarrow\infty} \Mb^{\gf_n}_{L,\la}(t)
\end{equation}
called the \textbf{stable fermionic formula}. Moreover, there are
precisely four distinct families $M^\se_{L,\la}(t)$ of such
polynomials, corresponding to the partitions $\se\in
\{\nn,\cell\,,\hdom\,,\vdom\,\}$ of size at most two.
\end{lemma}

See Figure \ref{fig:rootsys}. Note that two affine families have
the same associated partition $\se$ if and only if the
neighborhoods of the zero nodes of their Dynkin diagrams are
isomorphic.

\begin{remark}
The $M$ formula is used in the statement of Lemma \ref{lem:afflimit}
because it is well-defined for all $L$, whereas the $X$ formula
depends on the conjectural existence of KR crystals. The polynomial
$\Mb^{\gf_n}_{L,\la}(t)$ is equal to $M(W^L,\la,t^{-1})$ in the
notation of \cite{HKOTT01}. In the cases that the $X$ formula is
defined, it is not hard to see that the $X$ formula also has a large
rank limit. We call this polynomial $\Xb^\se_{L,\la}(t)$ the
\textbf{stable one-dimensional sum}.
\end{remark}

\subsection{The ubiquity of type $A$ in large rank}
Originally the $K$ formula was found by experimentation involving
graded analogues of certain creation operators for universal
characters of classical type \cite{SZ04}. The purpose of this
subsection is to give a representation-theoretic motivation for the
$K$ formula. It is the philosophy of M. Kleber, that the characters
of KR modules $W_s^{(r)}$ with $r$ nonspin, should behave like Schur
functions (type $A$ characters) \cite{Kl02}. One strong piece of
supporting evidence is that the KR characters conjecturally satisfy
a system of relations called the $Q$-system
\cite{HKOTT01,HKOTY99,KR90}. The $Q$-system relates KR characters
with others whose indices $r$ are nearby in the Dynkin diagram. In
the large rank limit the $Q$-system only involves the part of the
Dynkin diagram that looks like one of type $A$.

The implications of Kleber's idea for the $X=M$ conjecture,
can be made precise by the following two conjectures.

Given any dominant $\gb$-weight $\la$ involving no spin nodes, let
$W^\la$ be Chari and Pressley's minimal affinization \cite{CP95}. It
is a finite-dimensional $U'_q(\gf)$-module that need not have a
crystal base. In type $A^{(1)}$, $W^\la$ is isomorphic to the
$U_q(\gb)$-irreducible $V^\la$. In any nonexceptional type, the
character of the minimal affinization $W^\la$ should behave like the
Schur function $s_\la$.

\begin{conjecture} \label{conj:Atensor}
Fix $L=(L^{(r)}_s)$. For any nonexceptional affine algebra
of sufficiently large rank,
\begin{equation*}
  W^L \cong \bigoplus_\tau \Xb_{L,\tau}^A(1) W^\tau
\end{equation*}
up to filtration, where $\Xb_{L,\tau}^A(t)$ is the type $A$ stable
one-dimensional sum.
\end{conjecture}
Due to the definition of $X$ as a tensor product multiplicity,
this holds for type $A$, at least on the level of characters.

For $\se\in \theset$, let $\Par^\se$ denote the set of partitions
whose Ferrers diagrams may be tiled by the shape $\se$. Explicitly,
$\Par^\nn$ is the singleton set consisting of the empty partition,
$\Par^\cell$ is the set $\Par$ of all partitions, $\Par^\hdom$ is
the set of partitions with even rows, and $\Par^\vdom$ the set of
partitions with even columns.

\begin{conjecture} \label{conj:restrict} \cite{CK00} Let $\se\in\theset$ and
$\{\gf_n\}$ a nonexceptional family of affine algebras associated
with $\se$ as in Figure \ref{fig:rootsys}. Then for $n$ sufficiently
large and partitions representing dominant weights as above, one has
\begin{equation*}
  W^\la \cong \bigoplus_\tau \left(\sum_{\mu\in \Par^\se} c^\tau_{\la\mu}\right) V^\la
\end{equation*}
as $U_q(\gb)$-modules, where $c^\tau_{\la\mu}$ is the
Littlewood-Richardson (LR) coefficient or stable type A tensor
multiplicity (see subsection \ref{ss:LR}).
\end{conjecture}

Suppose $\la$ is the rectangular partition having $r$ rows and $s$
columns; this corresponds to the weight $s\omega_r$. It is expected
that the KR module $W^{(r)}_s$ is the minimal affinization
$W^{s\omega_r}$. In this case the above decomposition into
$U_q(\gb)$-modules, agrees with that which is prescribed by
\cite{HKOTT01,HKOTY99} for the KR module $W^{(r)}_s$.

Using the definition of $X$ as a branching tensor multiplicity,
Conjectures \ref{conj:Atensor} and \ref{conj:restrict} imply that
\begin{equation*}
  \Xb_{L,\la}^\se(1) = \sum_\tau \Xb^A_{L,\tau}(1) \sum_{\mu\in \Par^\se} c^\tau_{\la\mu}.
\end{equation*}
The $K$ formula may be obtained by putting a grading parameter into
this formula.

\subsection{The $K$ formula}
Define the polynomials $\Kb^\se_{L,\la}(t)\in\Z_{\ge0}[t]$ by
\cite{SZ04}
\begin{equation} \label{eq:K} \Kb^\se_{L,\la}(t) =
t^{|L|-|\la|} \sum_\tau \Xb^A_{L,\tau}(t^2) \sum_{\mu\in \Par^\se}
c^\tau_{\la\mu}
\end{equation}
where $|\la|=\sum_i \la_i$ and $|L|=\sum_{r,s} rs L^{(r)}_s$.

\begin{conjecture} \label{conj:create} \cite{SZ04} ($\mathbf{X=M=K}$)
For all $\se\in\{\nn,\,\cell\,,\,\vdom\,,\,\hdom\,\}$,
\begin{equation}
\label{eq:affch} \Kb^\se_{L,\la}(t) = \Mb^\se_{L,\la}(t^2)
\end{equation}
\end{conjecture}

The conjecture was previously known for $\se=\nn$ by combining the
results of a number of papers; see \cite{SZ04} for an explanation.
Our main result is:

\begin{theorem} \label{th:KX}
$\Xb^\se_{L,\la}(t^2)=\Kb^\se_{L,\la}(t)$ for $\se\in\onetwoset$ and
for $L$ such that $L_s^{(r)}=0$ for $r>1$, that is, when all tensor
factors have the form $B^{1,s}$.
\end{theorem}

Combining this with \cite{ScS04} we have $X=M=K$ for tensor products
of $B^{1,s}$ for large rank in types $\se\in\ourset$.

\begin{remark} \label{rem:xi}
Let
\begin{equation} \label{eq:xi}
\xi = \begin{cases} 2 & \text{if $\se=\cell\,$} \\
1 & \text{otherwise.}
\end{cases}
\end{equation}
Our coenergy function $\Db$ has value equal to $-\frac{1}{\xi}$
times the energy function $D$ used in \cite{HKOTT01}. Our $K$
formula in \cite{SZ04} uses energy instead of coenergy. The notation
$R$ in \cite{SZ04} refers to a list of rectangular partitions whose
multiplicities are given by $L$. The precise relationship is
\begin{equation} \label{eq:Kco}
\Kb^\se_{L,\la}(t) = t^{||L||+|L|-|\la|} K^\se_{\la,R}(t)
\end{equation}
where $||L||=\sum_{i,j\ge1} \binom{r_{ij}}{2}$ and
$r_{ij}=\sum_{p\ge i}\sum_{q\ge j} L_q^{(p)}$.
\end{remark}

\subsection{Combinatorial proof sketch} \label{ss:sketch}
By definition the $X$-formula has the form
\begin{equation} \label{eq:X}
  \Xb^\se_{L,\la}(t) = \sum_{b\in P_\se(L,\la)} t^{\Db(b)}
\end{equation}
where $P_\se(L,\la)$ is a finite set and
$\Db:P_\se(L,\la)\rightarrow \frac{1}{\xi}\Z_{\ge0}$ is a function.
They are described explicitly later.

Let $\LR(\tau;\la,\mu)$ be an appropriate set of cardinality
$c^\tau_{\la\mu}$; see section \ref{ss:LR}.

To prove Theorem \ref{th:KX} it suffices to exhibit a bijection
\begin{equation}\label{eq:bij}
\begin{split}
  \eta_L: P_\se(L,\la) &\rightarrow \bigcup_\tau \bigcup_{\mu\in \Par^\se}
  P_\nn(L,\tau) \times \LR(\tau;\la,\mu) \\
b\qquad &\mapsto \qquad(c,Z)
\end{split}
\end{equation}
such that
\begin{equation} \label{eq:bijstat}
  2\, \Db_\se(b) = |L| - |\la| + 2\, \Db_\nn(c).
\end{equation}

\subsection{Outline}
Sections \ref{sec:crystals} and \ref{sec:A} review affine crystal
theory and in particular that of type $A$. Section \ref{sec:VXR}
defines the map $\eta_L$ using the virtual crystal construction of
\cite{OSS03} and the combinatorial $R$-matrices. The virtual crystal
theory guarantees that this map respects the grading, provided it is
well-defined. The well-definedness and surjectivity of the map
require some work. Crystal embeddings called right-splitting are
used to reduce to the case of tensor powers of the vector KR crystal
$B^{1,1}$. In this case bijectivity and well-definedness are
established by considering in section \ref{sec:DDF} a completely
different realization of the bijection which we call the DDF map,
since its main ingredient is a bijection due to Delest, Dulucq, and
Favreau \cite{DDF88}. The DDF map has the advantage of being
obviously bijective, but it is not at all clear that it should
respect the grading. The proof is completed by observing inductively
that the VXR and DDF maps agree. In section \ref{sec:Roby} a more
intrinsic definition of the DDF map is given, extending the
interpretation of the original DDF bijection due to Roby \cite{R95}.

Section \ref{sec:end} contains a brief discussion about the case not
addressed in this article, namely, $\se=\vdom\,$.

While this manuscript was being completed, a proof appeared for the
ungraded $t=1$ case (with signed binomials) of the $X=M$ conjecture
\cite{Her05}. We were also informed by Cedric Lecouvey of his paper
\cite{Le04a}, which proves the special case of $X=K$ (Theorem
\ref{th:KX}) for tensor powers of $B^{1,1}$ in type $\hdom\,$. He
also proved in \cite{Le04} that for types $\hdom\,,\vdom\,$ and for
tensor products with factors of the form $B^{1,s}$, the $K$
polynomials are the affine Kazhdan-Lusztig polynomials of types $C$
and $D$ respectively, given by Lusztig's q-analogue of weight
multiplicity \cite{Lu83}. Combined with Theorem \ref{th:KX}, the
one-dimensional sums $X$ in type $\hdom$ given by tensor products of
crystals of the form $B^{1,s}$, are affine KL polynomials of type
$C$. Conjecturally the similar one-dimensional sums $X$ of type
$\vdom$ are affine KL polynomials of type $D$.

Many thanks to Cedric Lecouvey, Masato Okado, Anne Schilling, and
Mike Zabrocki for fruitful discussions and collaborations related to
this project.

\section{Review of affine crystal graphs}
\label{sec:crystals}

The definitions for general crystal graphs follow Kashiwara
\cite{Kas95}. The fundamental concepts for the study of crystal
graphs of finite-dimensional modules over affine algebras were
established in \cite{KMN92}.

\subsection{Affine algebras} Let $\gf\supset \gf'\supset \gb$
be an affine Kac-Moody algebra \cite{Kac90} with Dynkin vertex set
$I=\{0,1,\dotsc,n\}$, its derived subalgebra, and the simple Lie
subalgebra with Dynkin vertex set $J=I\backslash\{0\}$. Write
$U_q(\gf)\supset U'_q(\gf)\supset U_q(\gb)$ for the corresponding
quantized universal enveloping algebras \cite{KMN92}. For $\gf$ let
$\{\alpha_i\mid i\in I\}$ be the simple roots, $\{h_i\mid i\in I\}$
the simple coroots, $\{\Lambda_i\mid i\in I\}$ the fundamental
weights and $\delta$ the null root. Let
$P=\Z\delta\oplus\bigoplus_{i\in I} \Z\La_i$ and
$\Pb=\bigoplus_{i\in J} \Z\omega_i$ be the weight lattices for $\gf$
and $\gb$ respectively, where $\la\mapsto\bar{\la}$ denotes the
natural projection $P\rightarrow\Pb$ and $\{\omega_i=\bar{\La}_i\mid
i\in J\}$ are the fundamental weights for $\gb$.

\subsection{Affine crystal graphs}
Suppose $M$ is a finite-dimensional $U'_q(\gf)$-module with a
crystal base $B$. Then $B$ has the structure of a directed graph
with vertex set $B$ and directed edges colored by the set $I$, such
that:
\begin{enumerate}
\item If all edges are removed except those colored $i$, then the
connected components are finite directed paths called the
$i$-strings of $B$. Given $b\in B$, define $f_i(b)$ (resp. $e_i(b)$)
to be the vertex following (resp. preceding) $b$ in its $i$-string;
if there is no such vertex, declare the result to be the special
symbol $\emptyset$, which represents the zero element in the module
$M$. Define $\vp_i(b)$ (resp. $\ve_i(b)$) to be the number of arrows
from $b$ to the end (resp. beginning) of its $i$-string.
\item There is a weight function $\wt:B\rightarrow \Pb$ such that
for all $i\in I$,
\begin{equation*}
\begin{split}
\wt(f_i(b))&=\wt(b)-\overline{\alpha}_i \qquad\text{if $f_i(b)\not=\emptyset$}\\
\vp_i(b)-\ve_i(b) &= \inner{h_i}{\wt(b)}
\end{split}
\end{equation*}
\end{enumerate}
We shall call a colored directed graph $B$ with the above properties
a $U'_q(\gf)$-crystal.

A morphism $g:B\rightarrow B'$ of $U'_q(\gf)$-crystals is a map
$B\cup\{\emptyset\}\rightarrow B'\cup\{\emptyset\}$ such that for
any $b\in B$, either $g(b)\in B'$ or $g(b)=\emptyset$, and
$g(f_i(b))=f_i(g(b))$ and $g(e_i(b))=e_i(g(b))$ for all $b\in B$ and
$i\in I$ where by definition $g(\emptyset)=\emptyset$. An
isomorphism of crystals is a morphism of crystals which is a
bijection whose inverse is also a morphism of crystals.

The direct sum of crystals is the disjoint union.

\subsection{Crystal graphs for simple Lie algebras}
If $M$ is a finite-dimensional $U_q(\gb)$-module with crystal base
$B$ then its directed edges are colored by the set
$J=I\backslash\{0\}$, and for $i\in J$ it satisfies the properties
stated previously for affine crystals. In this case we call $B$ a
$U_q(\gb)$-crystal.

$b\in B$ is a \textbf{classical highest weight vector} if
$\ve_i(b)=0$ for all $i\in J$.

For $\la\in \Pb^+$ let $B(\la)$ be the crystal graph of the
irreducible highest weight $U_q(\gb)$-module of highest weight
$\la$. It has a unique classical highest weight vector and this
vector has weight $\la$. For $\gb$ of classical type $B(\la)$ was
given explicitly in \cite{KN94}.

Every finite-dimensional $U_q(\gb)$-module has a crystal graph which
is isomorphic to a direct sum of crystal graphs of the form
$B(\la)$.

A \textbf{classical component} of a $U'_q(\gf)$-crystal $B$ is a
connected component of the $U_q(\gb)$-crystal obtained from $B$ by
removing all the $0$-arrows.

\subsection{Tensor products}
\begin{remark}
We use the opposite of Kashiwara's tensor product convention to be
consistent with the $U_q(A_n)$-crystal graph combinatorics of Young
tableaux.
\end{remark}

If $B_1,B_2$ are crystal bases of the $U'_q(\gf)$-modules $M_1,M_2$
then the tensor product $M_2 \otimes M_1$ has a crystal base denoted
$B_2\otimes B_1$. Its vertex set is just the cartesian product $B_2
\times B_1$. Its edges are given by
\begin{equation} \label{eq:ftensor}
\begin{split}
  f_i(b_2\otimes b_1) &= \begin{cases}
  f_i(b_2)\otimes b_1 &\text{if $\ve_i(b_2)\ge\vp_i(b_1)$} \\
  b_2 \otimes f_i(b_1) &\text{otherwise.}
\end{cases} \\
  e_i(b_2\otimes b_1) &= \begin{cases}
  e_i(b_2)\otimes b_1 &\text{if $\ve_i(b_2)>\vp_i(b_1)$} \\
  b_2 \otimes e_i(b_1) &\text{otherwise,}
\end{cases}
\end{split}
\end{equation}
where the result is declared to be $\emptyset$ if either of its
tensor factors are. The weight function on $B_2\otimes B_1$ is
$\wt(b_2\otimes b_1)=\wt(b_2)+\wt(b_1)$.

The tensor product operation on crystals is associative.

\subsection{Contragredient dual} \label{ss:dual}
Given the crystal base $B$ of a $U'_q(\gf)$-module $M$, its
contragredient dual module $M^\vee$ has a crystal base $B^\vee$,
which has vertices written $b^\vee$ for $b\in B$, with
$\wt(b^\vee)=-\wt(b)$ and with all arrows reversed \cite{Kas95}:
\begin{equation} \label{eq:veeraise}
  f_i(b^\vee) = e_i(b)^\vee\qquad\text{for all $i\in I$ and $b\in
  B$.}
\end{equation}
There is a natural $U'_q(\gf)$-crystal isomorphism
\begin{equation} \label{eq:veetensor}
\begin{split}
 (B_L\otimes \dotsm\otimes B_1)^\vee &\cong B_1^\vee \otimes\dotsm\otimes
 B_L^\vee \\
 (b_L\otimes\dotsm\otimes b_1)^\vee &\mapsto b_1^\vee \otimes\dotsm\otimes
 b_L^\vee
\end{split}
\end{equation}

\subsection{Representative affine families}
\begin{figure}
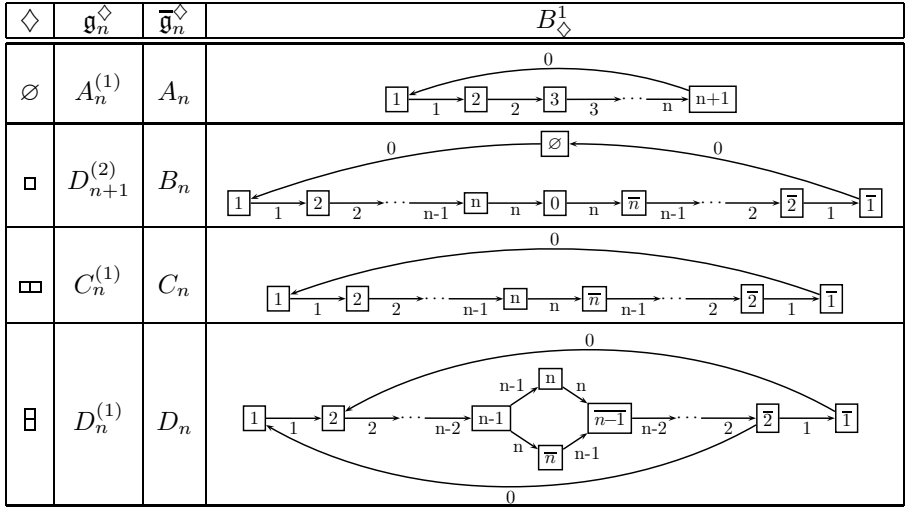

\begin{tabular}{|c|c|c|c|}
\hline
$\se$ & $\gf_n^\se$ & $\gb_n^\se$ & $B_\se^1$ \\
\hline \hline
%
$\nn$  & $A_n^{(1)}$ & $A_n$ & \raisebox{-3mm}{%
\scalebox{.7}{%
\psset{unit=1.5cm,labelsep=.5mm}
\pspicture(0.75,0.7)(5.25,1.7)%
\rput(1,1){\rnode{A1}{\psframebox{1}}}%
\rput(2,1){\rnode{A2}{\psframebox{2}}}%
\rput(3,1){\rnode{A3}{\psframebox{3}}}%
\rput(4,1){\rnode{A4}{$\dotsm$}}%
\rput(5,1){\rnode{A5}{\psframebox{n+1}}}%
\ncline{->}{A1}{A2}%
\nbput{1}%
\ncline{->}{A2}{A3}%
\nbput{2}%
\ncline{->}{A3}{A4}%
\nbput{3}%
\ncline{->}{A4}{A5}%
\nbput{n}%
\nccurve[angleA=160,angleB=20]{->}{A5}{A1}%
\nbput{0}%
\endpspicture
}}%
\\ \hline
%
$\cell$ & $D_{n+1}^{(2)}$ & $B_n$ & \raisebox{-5mm}{%
\scalebox{.7}{%
\psset{unit=1.5cm,labelsep=.5mm}
\pspicture(0.75,0.7)(9.25,2)%
\rput(1,1){\rnode{A1}{\psframebox{1}}}%
\rput(2,1){\rnode{A2}{\psframebox{2}}}%
\rput(3,1){\rnode{A3}{$\dotsm$}}%
\rput(4,1){\rnode{A4}{\psframebox{n}}}%
\rput(5,1){\rnode{A5}{\psframebox{0}}}%
\rput(6,1){\rnode{A6}{\psframebox{$\overline{n}$}}}%
\rput(7,1){\rnode{A7}{$\dotsm$}}%
\rput(8,1){\rnode{A8}{\psframebox{$\overline{2}$}}}%
\rput(9,1){\rnode{A9}{\psframebox{$\overline{1}$}}}%
\rput(5,1.75){\rnode{B}{\psframebox{$\varnothing$}}}%
\ncline{->}{A1}{A2}%
\nbput{1}%
\ncline{->}{A2}{A3}%
\nbput{2}%
\ncline{->}{A3}{A4}%
\nbput{n-1}%
\ncline{->}{A4}{A5}%
\nbput{n}%
\ncline{->}{A5}{A6}%
\nbput{n}%
\ncline{->}{A6}{A7}%
\nbput{n-1}%
\ncline{->}{A7}{A8}%
\nbput{2}%
\ncline{->}{A8}{A9}%
\nbput{1}%
\nccurve[angleA=160,angleB=0]{->}{A9}{B}%
\nbput{0}%
\nccurve[angleA=180,angleB=20]{->}{B}{A1}%
\nbput{0}%
\endpspicture
}}%

\\ \hline
%
$\hdom$ & $C_n^{(1)}$ & $C_n$ & \raisebox{-4mm}{%
\scalebox{.7}{%
\psset{unit=1.5cm,labelsep=.5mm}
\pspicture(0.75,0.7)(8.25,1.9)%
\rput(1,1){\rnode{A1}{\psframebox{1}}}%
\rput(2,1){\rnode{A2}{\psframebox{2}}}%
\rput(3,1){\rnode{A3}{$\dotsm$}}%
\rput(4,1){\rnode{A4}{\psframebox{n}}}%
\rput(5,1){\rnode{A5}{\psframebox{$\overline{n}$}}}%
\rput(6,1){\rnode{A6}{$\dotsm$}}%
\rput(7,1){\rnode{A7}{\psframebox{$\overline{2}$}}}%
\rput(8,1){\rnode{A8}{\psframebox{$\overline{1}$}}}%
\ncline{->}{A1}{A2}%
\nbput{1}%
\ncline{->}{A2}{A3}%
\nbput{2}%
\ncline{->}{A3}{A4}%
\nbput{n-1}%
\ncline{->}{A4}{A5}%
\nbput{n}%
\ncline{->}{A5}{A6}%
\nbput{n-1}%
\ncline{->}{A6}{A7}%
\nbput{2}%
\ncline{->}{A7}{A8}%
\nbput{1}%
\nccurve[angleA=160,angleB=20]{->}{A8}{A1}%
\nbput{0}%
\endpspicture
}}%
\\ \hline

$\vdom$ & $D_n^{(1)}$ & $D_n$ & \raisebox{-1cm}{%
\scalebox{.7}{%
\psset{unit=1.5cm,labelsep=.5mm}
\pspicture(-0.25,-1.1)(7.75,1.2)%
\rput(0,0){\rnode{A1}{\psframebox{1}}}%
\rput(1,0){\rnode{A2}{\psframebox{2}}}%
\rput(2,0){\rnode{A3}{$\dotsm$}}%
\rput(3,0){\rnode{A4}{\psframebox{n-1}}}%
\rput(3.75,.5){\rnode{B1}{\psframebox{n}}}%
\rput(3.75,-.5){\rnode{B2}{\psframebox{$\overline{n}$}}}%
\rput(4.5,0){\rnode{A5}{\psframebox{$\overline{n\!\!-\!\!1}$}}}%
\rput(5.5,0){\rnode{A6}{$\dotsm$}}%
\rput(6.5,0){\rnode{A7}{\psframebox{$\overline{2}$}}}%
\rput(7.5,0){\rnode{A8}{\psframebox{$\overline{1}$}}}%
\ncline{->}{A1}{A2}%
\nbput{1}%
\ncline{->}{A2}{A3}%
\nbput{2}%
\ncline{->}{A3}{A4}%
\nbput{n-2}%
\ncline{->}{A4}{B1}%
\naput{n-1}%
\ncline{->}{A4}{B2}%
\nbput{n}%
\ncline{->}{B1}{A5}%
\naput{n}%
\ncline{->}{B2}{A5}%
\nbput{n-1}%
\ncline{->}{A5}{A6}%
\nbput{n-2}%
\ncline{->}{A6}{A7}%
\nbput{2}%
\ncline{->}{A7}{A8}%
\nbput{1}%
\nccurve[angleA=150,angleB=30]{->}{A8}{A2}%
\nbput{0}%
\nccurve[angleA=-150,angleB=-30]{->}{A7}{A1}%
\naput{0}%
\endpspicture
}
}%
\\ \hline
\end{tabular}\vspace{4mm}
\caption{\label{fig:crystals} Representative affine families and
$B_\se^1$}
\end{figure}

For each partition $\se\in \{\nn,\cell\,,\hdom\,,\vdom\,\}$, we fix
a representative family of affine root systems. An affine algebra in
this family and its distinguished simple Lie subalgebra shall be
denoted $\gf_n^\se$ and $\gb_n^\se$ respectively. The subscript $n$
indicates the rank of $\gb_n^\se$. See Figure \ref{fig:crystals}.

\subsection{Categories of affine crystals}
\label{ss:assume} We use two categories of affine crystal graphs:
the category $\CC(A_n^{(1)})$ of tensor products of arbitrary KR
crystals of type $A_n^{(1)}$ and the category $\CC^1(\gf^\se_n)$ of
tensor products of KR crystals of the form $B_\se^s=B^{1,s}$ for
$\gf_n^\se$.

Let $\CC=\CC(A_n^{(1)})$ or $\CC=\CC^1(\gf_n^\se)$. $\CC$ satisfies
the following properties \cite{AK97,KKM94}:
\begin{enumerate}
\item Every $B\in \CC$ is the crystal base of an irreducible
finite-dimensional $U'_q(\gf)$-module.
\item $\CC$ is closed under tensor product.
\item For every $B\in\CC$ there is a $\la_0\in \Pb^+$ such that there is a
unique element $u(B)\in B$ of weight $\la_0$, and such that all
other weights occurring in $B$ are in the convex hull of
$\overline{W} \la_0$ where $\overline{W}$ is the Weyl group of
$\gb$. We call $u(B)$ the \textbf{leading vector} of $B$ and $\la_0$
the \textbf{leading weight} of $B$.
\end{enumerate}
The leading weight of $B^{r,s}$ is $s\omega_r$.

\subsection{KR crystals} \label{ss:KR}
We give the KR crystals in the categories $\CC$.

\subsubsection{Vector KR crystals} \label{sss:boxKR}
For $\gf=\gf_n^\se$, the KR crystal $B_\se^1$ is depicted in Figure
\ref{fig:crystals}. It is a $U'_q(\gf)$-crystal with
$U_q(\gb)$-decomposition given by
\begin{equation}
  B_\se^1 = \begin{cases}
  B(\omega_1)\oplus B(0) & \text{if $\se=\cell\,$} \\
  B(\omega_1) &\text{otherwise.}
  \end{cases}
\end{equation}
Removing zero arrows from $B_\se^1$ and in the case $\se=\cell\,$
removing the element $\nn\in B(0)$, one obtains the
$U_q(\gb)$-crystal graph $B(\omega_1)$.

Define a partial order on $B(\omega_1)$ by $x<y$ if there is a
directed path in $B(\omega_1)$ from $x$ to $y$.

The tensor product rule \eqref{eq:ftensor} gives the
$U_q(\gb)$-crystal structure on tensor powers of $B(\omega_1)$.

\subsubsection{Row KR crystals}
\label{sss:rowKR} More generally the KR crystal $B_\se^s$ is a
$U'_q(\gf)$-crystal with $U_q(\gb)$-decomposition given by
\cite{KKM94}
\begin{equation} \label{eq:Bdecomp}
  B^s_\se =
\begin{cases}
B(s\omega_1) & \text{if $\se\in\{\nn,\vdom\,\}$} \\
\bigoplus_{r=0}^s B((s-r)\omega_1) &\text{if $\se=\cell\,$} \\
\bigoplus_{r=0}^{\lfloor s/2\rfloor} B((s-2r)\omega_1) &\text{if
$\se=\hdom\,$.}
\end{cases}
\end{equation}
To give the $U_q(\gb)$-crystal structure on $B_\se^s$ it suffices to
define such a structure for its classical components, all of which
have the form $B(m\omega_1)$. As a set $B(m\omega_1)$ consists of
the weakly increasing sequences of length $m$ of elements in
$B(\omega_1)$, with the additional constraints that $0$ may only
occur once for $\se=\cell\,$ and that $n$ and $\bar{n}$ may not both
appear for $\se=\vdom\,$. The $U_q(\gb)$-crystal structure on
$B(m\omega_1)$ is given by observing that the following map is a
$U_q(\gb)$-crystal embedding
\begin{equation} \label{eq:rowtensor}
\begin{split}
B(m\omega_1)&\overset{\iota_m}{\longrightarrow} B(\omega_1)^{\otimes m} \\
b_m\dotsm b_1&\mapsto b_m\otimes\dotsm\otimes b_1
\end{split}
\end{equation}
where $b_j\in B(\omega_1)$.

To specify the full $U'_q(\gf)$-crystal structure on $B_\se^s$ one
must also know the $0$-arrows. This is given explicitly in
\cite{KKM94,OSS03b}. Rather than recalling this here we shall define
it indirectly using the virtual crystal construction
\cite{OSS03,OSS03b}. This is postponed until section \ref{ss:Vrow}.

\subsection{Combinatorial $R$-matrix} Given any $B_1,B_2\in \CC$
with underlying $U'_q(\gf)$-modules $M_1,M_2$, the universal
$R$-matrix induces a $U'_q(\gf)$-module isomorphism $M_2\otimes
M_1\rightarrow M_1\otimes M_2$. By assumption the crystal graphs of
these tensor products are connected. It follows that there is a
unique $U'_q(\gf)$-crystal isomorphism $R=R_{B_2,B_1}:B_2\otimes B_1
\rightarrow B_1\otimes B_2$ called the \textbf{combinatorial
$R$-matrix}. Since $R$ preserves weights, leading vectors must
correspond: $u(B_2\otimes B_1)=u(B_2)\otimes u(B_1)\mapsto
u(B_1\otimes B_2)=u(B_1)\otimes u(B_2)$. In particular uniqueness
implies that
\begin{equation} \label{eq:R=}
R_{B,B}=1_{B\otimes B}.
\end{equation}
For $\CC^1(\gf)$ for $\gf$ nonexceptional, the combinatorial
$R$-matrices are given by \cite{HKOT00,HKOT02}. For $\CC(A_n^{(1)})$
the combinatorial $R$-matrix is discussed in section \ref{ss:KRA}.

\subsection{Local coenergy} To agree with the convention in
\eqref{eq:affch} and \eqref{eq:xi}, let $\xi=1$ except if, in the
Dynkin diagram of $\gf$ there is a double arrow pointing toward the
$0$ node, in which case $\xi=2$. The only nonexceptional affine
families with $\xi=2$ are $D_{n+1}^{(2)}$ and $A_{2n}^{(2)}$; these
are also precisely the affine families associated with $\cell\,$ in
Figure \ref{fig:rootsys}. By \cite{KMN92}, for all $B_1,B_2\in \CC$
there is a unique function $\Hb=\Hb_{B_2\otimes B_1}:B_2\otimes
B_1\rightarrow \frac{1}{\xi}\Z_{\ge0}$ called the \textbf{local
coenergy function}, such that
\begin{enumerate}
\item $\Hb$ is constant on classical components of $B_2\otimes B_1$.
\item By \eqref{eq:ftensor}, for $b_2\otimes b_1\in B_2\otimes
B_1$, $e_0(b_2\otimes b_1)$ equals either $e_0(b_2)\otimes b_1$ or
$b_2\otimes e_0(b_1)$. The same is true for $e_0(b_1'\otimes b_2')$
where $R_{B_2,B_1}(b_2\otimes b_1)=b_1'\otimes b_2'$. Say that $LL$
(resp. RR) holds if $e_0$ acts on the left (resp. right) tensor
factor for both $b_2\otimes b_1$ and $b_1'\otimes b_2'$. Then
\begin{equation} \label{eq:H}
  \Hb(e_0(b_2\otimes b_1)) = \Hb(b_2\otimes b_1)+ \dfrac{1}{\xi}
  \cdot
\begin{cases}
1 & \text{if $LL$} \\
-1 & \text{if $RR$} \\
0 &\text{otherwise.}
\end{cases}
\end{equation}
\item $\Hb(u(B_2)\otimes u(B_1))=0$.
\end{enumerate}
By definition
\begin{equation} \label{eq:HR}
\Hb_{B_1\otimes B_2} = \Hb_{B_2\otimes B_1} \circ R_{B_1,B_2}.
\end{equation}
If $B_1=B_2=B$ then by \eqref{eq:R=}, \eqref{eq:H} simplifies to
\begin{equation}\label{eq:H=}
  \Hb(e_0(b_2\otimes b_1))=\Hb(b_2\otimes b_1)+
\dfrac{1}{\xi}
  \cdot\begin{cases}
1 &\text{if $e_0(b_2\otimes b_1)=e_0(b_2)\otimes b_1$} \\
-1 &\text{otherwise.}
\end{cases}
\end{equation}
For $\CC^1(\gf)$ the local coenergy function is given explicitly in
\cite{HKOT00,HKOT02} and for $\CC(A_n^{(1)})$ it is given in
\cite{S02}.

\begin{example} \label{ex:Hbox}
Using the partial order on $B(\omega_1)$ given in section
\ref{sss:boxKR} we have
\begin{equation} \label{eq:HB11}
\begin{split}
  \Hb_{B_\nn\otimes B_\nn}(x \otimes y) &=
  \begin{cases}
    0 & \text{if $x\le y$} \\
    1 & \text{otherwise}
  \end{cases} \\
  \Hb_{B_\hdom\,\otimes B_\hdom}(x \otimes y) &=
  \begin{cases}
    0 & \text{if $x\le y$} \\
    1 & \text{otherwise}
  \end{cases} \\
  \Hb_{B_\vdom\,\otimes B_\vdom}(x \otimes y) &=
  \begin{cases}
    2& \text{if $(x,y)=(\bar{1},1)$} \\
    0& \text{if $x\le y$} \\
    1& \text{otherwise}
  \end{cases} \\
  \Hb_{B_\cell\, \otimes B_\cell}(x\otimes y) &=
  \begin{cases}
     1&\text{if $x>y$ or $(x,y)=(\nn,\nn)$ or $(x,y)=(0,0)$} \\
     \frac{1}{2}&\text{if exactly one of $x$ or $y$ is $\nn$} \\
     0&\text{otherwise.}
  \end{cases}
\end{split}
\end{equation}
\end{example}

\subsection{Yang-Baxter equation} \label{ss:YB}
This section follows \cite{KMN92}.  Given $B_1,B_2,B_3\in \CC$,
there is a unique isomorphism $B_3 \otimes B_2 \otimes
B_1\rightarrow B_1\otimes B_2\otimes B_3 $ induced by the algebraic
universal $R$-matrix. By uniqueness the combinatorial $R$-matrices
satisfy the Yang-Baxter equation
\begin{equation} \label{eq:YB}
  R_{B_3,B_2} \circ R_{B_3,B_1}\circ R_{B_2,B_1} =
  R_{B_2,B_1} \circ R_{B_3,B_1} \circ R_{B_3,B_2}.
\end{equation}
Using sloppier notation, let $R_i$ denote the $R$-matrix which acts
in the $i$-th and $(i+1)$-st positions from the right in a tensor
product. Then we may rewrite the Yang-Baxter equation as the braid
relation $R_1 R_2 R_1=R_2 R_1 R_2$.

Again by uniqueness one has
\begin{equation} \label{eq:R21}
\begin{split}
  R_{(B_1 \otimes B_2) \otimes B_3} &= (R_{B_1\otimes B_3} \otimes
  1_{B_2}) \circ (1_{B_1} \otimes R_{B_2\otimes B_3}) \\
  R_{B_1 \otimes (B_2 \otimes B_3)} &=
  (1_{B_2} \otimes R_{B_1\otimes B_3}) \circ
  (R_{B_1 \otimes B_2} \otimes 1_{B_3})
\end{split}
\end{equation}

Let $H_i$ denote the value of the local coenergy function, evaluated
at the $i$-th and $(i+1)$-st positions from the right in a tensor
product. Then one has
\begin{equation} \label{eq:HYB}
\begin{split}
  H_1 + H_2 R_1 &= H_1 R_2 + H_2 R_1 R_2 \\
  H_2 + H_1 R_2 &= H_2 R_1 + H_1 R_2 R_1
\end{split}
\end{equation}

The following consequence of \eqref{eq:YB} and \eqref{eq:HYB} is
well-known.

\begin{prop} \label{pp:Hbraided}
Let $b\in B=B_L\otimes\dotsm\otimes B_1$ and $b'\in
B'=B_M'\otimes\dotsm \otimes B_1'$.
\begin{enumerate}
\item $R_{B \otimes B'}(b\otimes b')$ may be computed by a composition of
$R$-matrices of the form $R_{B_i \otimes B_j'}$ to shuffle the
element of $B_i$ to the right past that in $B_j'$.
\item $\Hb_{B \otimes B'}(b\otimes b')$ is the sum of
the values $\Hb_{B_i \otimes B_j'}$ evaluated at the pairs of
elements in $B_i \otimes B_j'$ that must be switched with an
$R$-matrix $R_{B_i \otimes B_j'}$ in the computation of $R_{B
\otimes B'}(b\otimes b')$.
\end{enumerate}
\end{prop}

\subsection{Intrinsic coenergy}
\label{ss:D} We give a construction introduced in \cite{HKOTT01} as
described in \cite{OSS03} for $\CC$ as in section \ref{ss:assume}.
For $B\in \CC$, we shall define the intrinsic coenergy function
$\Db_B:B\rightarrow \frac{1}{\xi}\Z_{\ge0}$. It satisfies the
following rules.
\begin{enumerate}
\item $\Db_B$ is constant on classical components.
\item $\Db_B(e_0(b))\le 1/\xi+\Db_B(b)$ for all $b\in B$.
\item $\Db_B(u(B))=0$.
\end{enumerate}

\begin{example} \label{ex:KRD}
For KR crystals $B^{r,s}$ the intrinsic coenergy functions are
prescribed in \cite{HKOTT01,HKOTY99}. For $\se\in\onetwoset$,
$\gf_n^\se$ as in Figure \ref{fig:crystals}, and $B_\se^s$ with
classical components as in \eqref{eq:Bdecomp}, the $r$-th tensor
factor has intrinsic coenergy $r/\xi$ with $\xi$ as in
\eqref{eq:xi}. In particular,
\begin{equation} \label{eq:DB11}
  \Db_{B_\se^1}(x) = \begin{cases}
   \frac{1}{2} &\text{if $\se=\cell\,$ and $x=\nn$} \\
   0 &\text{otherwise.}
  \end{cases}
\end{equation}
\end{example}

Suppose the intrinsic coenergy functions $D_1$ and $D_2$ of
$B_1,B_2\in\CC$ have already been defined. Define the intrinsic
coenergy function for $B_2\otimes B_1$ as follows. Let $\Db^1_{B_i}$
denote the function $\Db_{B_i}$ applied to the rightmost tensor
factor. Then
\begin{equation} \label{eq:D2}
  \Db_{B_2\otimes B_1} =
  \Hb_{B_2\otimes B_1} +
  \Db_{B_1}^1 +
  \Db_{B_2}^1 \circ R_{B_2,B_1}.
\end{equation}

\begin{prop} \label{pp:Dtensorcat} \cite{OSS03} The above
construction of $\Db$ is associative.
\end{prop}

Using this construction we may define $\Db_B$ for all $B\in \CC$. By
Proposition \ref{pp:Dtensorcat} the function $\Db_B$ doesn't depend
on the way that $B$ is built up from two-fold tensor products.

For $B_1,\dotsc,B_L\in\CC$, the above construction yields the
formula  \cite{HKOTT01,NY97,OSS03}
\begin{equation} \label{eq:Dn}
  \Db_{B_L\otimes\dotsm \otimes B_1} =
\sum_{1\le i<j\le L} \Hb_i R_{i+1} R_{i+2} \dotsm R_{j-1} +
\sum_{j=1}^L \Db_{B_j}^1 R_1 R_2 \dotsm R_{j-1}.
\end{equation}
If $B_j=B$ for all $j$ then by \eqref{eq:R=} this reduces to
\begin{equation} \label{eq:Dn=}
  \Db_{B^{\otimes L}} =
\sum_{i=1}^{L-1} (L-i)\, \Hb_i + L \,\Db_B^1.
\end{equation}

\begin{prop} \label{pp:DR} \cite{OSS03} Suppose
$B$ and $B'$ are any two tensor products with the same collection of
tensor factors $B_1,\dotsc,B_L\in\CC$ in some order. Let
$R:B\rightarrow B'$ be the unique affine crystal isomorphism. Then
\begin{equation} \label{eq:DR}
  \Db_B = \Db_{B'}\circ R
\end{equation}
\end{prop}

\subsection{The $X$ formula}
Let $B\in\CC^1(\gf_n^\se)$ and $\la\in\Pb^+$. Let $P(B)$ be the set
of classical highest weight vectors in $B$ and $P(B,\la)$ those of
highest weight $\la$. For type $\se$ we write $P_\se(B,\la)$. If
$B=B^{1\otimes L}$ we write $P_\se(1^L,\la)$.

\begin{example} \label{ex:standardpath}
Let $B=B_\se^{1\otimes L}$. Let $\gf$ be of type
$A_{n-1}^{(1)},D_{n+1}^{(2)},C_n^{(1)},D_n^{(1)}$ so that $\gb$ has
type $A_{n-1},B_n,C_n,D_n$ for types $\se=\theset$ respectively.

We use the identification of weights and elements of $\Z^n$ given in
section \ref{ss:limit}. Let $m_x(b)$ be the number of times the
symbol $x$ occurs in the word $b$. Define the weight function
$\wt:B_\se^{1\otimes L}\rightarrow\Z^n$ by
\begin{equation}
\begin{split}
  \wt(b) &=
  (m_1(b)-m_{\bar{1}}(b),m_2(b)-m_{\bar{2}}(b),\dotsc,m_n(b)-m_{\bar{n}}(b))
\end{split}
\end{equation}
For rank large with respect to $L$ and $\la$ it is easy to check
that $b\in P_\se(1^L,\la)$ if and only if $\la^{(i)}=\wt(b_i\dotsm
b_1)$ is a partition for all $i$, with $\la^{(L)}=\la$. In type $A$
such classical highest weight vectors are called \textbf{Yamanouchi
words}.
\end{example}

The \textbf{one-dimensional sum} is defined by
\begin{equation*}
  \Xb(B,\la)(t) = \sum_{b\in P(B,\la)} t^{\Db_B(b)}.
\end{equation*}

\begin{remark} \label{rem:Border}
By Proposition \ref{pp:DR} this depends only on the collection of
tensor factors in $B$ and not on their order. Let $L$ be the
multiplicity array for $B$, that is, $B$ contains $L_i^{(a)}$ tensor
factors equal to $B^{a,i}$. We write $\Xb(B,\la)(t)=\Xb_{L,\la}(t)$.
\end{remark}

\begin{example} \label{ex:X} We give examples of summands in
one-dimensional sums. Tensor symbols are dropped and subscripts are
added for convenience. In these examples $\Db$ is calculated using
\eqref{eq:Dn=}, \eqref{eq:HB11}, and \eqref{eq:DB11}.
\begin{equation*}
c=1_13_22_31_42_51_61_7\in P_\nn(1^7,(4,2,1)) \qquad
b=\bar{2}_11_2\nn_3\bar{1}_42_51_61_7 \in P_\cell\,(1^7,(2))
\end{equation*}
\begin{equation*}
\begin{split}
 \Db_\nn(c) &= 1\cdot0+2\cdot1+3\cdot1+4\cdot0+5\cdot1+6\cdot0+7\cdot 0  = 10 \\
 \Db_\cell\,(b) &= 1\cdot1+2\cdot\frac{1}{2}+3\cdot\frac{1}{2}+4\cdot1+5\cdot1+6\cdot0+7 \cdot 0=25/2.
\end{split}
\end{equation*}
\end{example}

\section{Type A crystal graphs}
\label{sec:A}

Let $N=2n$ for the remainder of the paper. In this section
$\gf=A_{N-1}^{(1)}$ and $\gb=A_{N-1}$. We specialize to this case
because of the virtual crystal construction in section
\ref{sec:VXR}.

\subsection{$U_q(A_{N-1})$-crystals}
\label{ss:Aclassical} The $U_q(A_{N-1})$-crystal structure on
$B(\omega_1)^{\otimes m}$ is given in section \ref{sss:boxKR}; it is
the set of words of length $m$ in the totally-ordered set
$\{1<2<\dotsm<N\}$.

The crystal graph $B(\la)$ of the highest weight
$U_q(A_{N-1})$-module of highest weight $\la$, is given by the set
of (semistandard) Young tableaux of shape $\la$ with entries in the
set $B(\omega_1)=\{1,2,\dotsc,N\}$ \cite{KN94}; see \cite{Ful97} for
the definition of a tableau and its row word.

Let $|\la|=m$ be the number of cells in the diagram of the partition
$\la$. The $U_q(A_{N-1})$-crystal structure on $B(\la)$ is given by
declaring that the map $B(\la)\rightarrow B(\omega_1)^{\otimes m}$
that sends a tableau to its row word, is an embedding of
$U_q(A_{N-1})$-crystals.

Let $u\in B(\omega_1)^{\otimes m}$. There is a unique dominant
$A_{N-1}$-weight $\la$ such that the component of $u$ in
$B(\omega_1)^{\otimes m}$ is isomorphic to $B(\la)$. Moreover the
isomorphism is unique. Denote by $P(u)\in B(\la)$ the image of $u$
under this isomorphism. It is well-known that $P(u)$ is Schensted's
$P$-tableau \cite{Sc61}.

\subsection{Knuth equivalence}
We shall identify a tableau with its row word. Knuth defined an
equivalence relation on words with letters in the totally-ordered
finite set $A$, which is generated by relations of the form
\begin{align*}
uxzyv&\Knuth uzxyv && \text{for $x<y\le z$} \\
uyxzv&\Knuth uyzxv && \text{for $x\le y<z$}
\end{align*}
where $u$ and $v$ are words with letters in $A$ and $x,y,z\in A$
\cite{Kn70}.

One may show that for each word $u$, there is a unique tableau of
partition shape that is Knuth-equivalent to $u$; we denote it by
$[u]$.

The following results are well-known.

\begin{prop} \label{pp:KnuthA} For any word $u$, $P(u)=[u]$.
\end{prop}

\begin{prop} \label{pp:Knuthres} Suppose $u\Knuth v$ where $u$ and $v$
are words with letters in a totally ordered set $A$ and $B$ is a
subinterval of $A$. Then $u|_B\Knuth v|_B$ where $u|_B$ is obtained
from $u$ by erasing all letters not in $B$.
\end{prop}

\subsection{Littlewood-Richardson Rule}
\label{ss:LR} The Littlewood-Richardson (LR) rule \cite{LR34} is a
combinatorial description of the structure constants of the ring of
symmetric functions with respect the basis of Schur functions. See
\cite{FG93} for an exposition of various combinatorial viewpoints on
this rule. It is well-known that $c^\alpha_{\beta\gamma}$ is the
multiplicity of the irreducible $U_q(A_{N-1})$-crystal $B(\alpha)$
in the tensor product $B(\beta)\otimes B(\gamma)$ for $N$ large with
respect to the partitions $\alpha,\beta,\gamma$.

Let $y$ be a word of length $m$ in the alphabet $\{1,2,\dotsc,N\}$,
that is, $y$ is in the $U_q(A_{N-1})$-crystal $B(\omega_1)^{\otimes
m}$. The word $y$ is \textbf{Yamanouchi} if it is a highest weight
vector. This is equivalent to saying that the weight of each right
factor of $y$ is a partition (and in particular $\wt(y)$ is a
partition). Say that a word $w$ is $\alpha/\beta$-Yamanouchi if for
any Yamanouchi word $y$ of weight $\beta$, $wy$ is Yamanouchi of
weight $\alpha$.

\begin{theorem} \label{th:LR} The Littlewood-Richardson
coefficient $c^\tau_{\la\mu}$ is the cardinality of any of the
following sets:
\begin{enumerate}
\item The tableaux of skew shape $\tau/\la$ that are Yamanouchi of
weight $\mu$ \cite{LR34}.
\item The tableaux of shape $\mu$ that are $\tau/\la$-Yamanouchi.
\item Given a fixed tableau $P$ of shape $\tau$,
the set of pairs of tableaux $(T,S)$ of respective shapes $\la$ and
$\mu$, such that $P=[TS]$ \cite{LS81}.
\end{enumerate}
\end{theorem}

\subsection{Type A KR crystals}
\label{ss:KRA} The KR $U'_q(A_{N-1}^{(1)})$-crystal $B^{r,s}$ is
isomorphic as a $U_q(A_{N-1})$-module to $B(s\omega_r)$
\cite{KMN92a}. By section \ref{ss:Aclassical} $B^{r,s}$ consists of
the semistandard tableaux with entries in $\{1,2,\dotsc,N\}$ whose 0
shape is the $r\times s$ rectangle. The $0$-arrows for $B^{r,s}$ are
given in \cite{KMN92a,S02}.

Since $B^{r,s} \otimes B^{r',s'}$ is multiplicity-free as a
$U_q(A_{N-1})$-crystal, the above combinatorial $R$-matrix
$R_{B^{r,s},B^{r',s'}}$ may be characterized as the unique
$U_q(A_{N-1})$-crystal isomorphism $B^{r,s}\otimes
B^{r',s'}\rightarrow B^{r',s'}\otimes B^{r,s}$ and can therefore be
computed using the Robinson-Schensted-Knuth correspondence
\cite{S02}.

We recall from \cite{S02} that for $b_1\otimes b_2\in B^{r,s}
\otimes B^{r',s'}$, the value of the local coenergy function
$\Hb(b_1 \otimes b_2)$ is the number of cells in the shape of the
tableau $[b_1 b_2]$ that are strictly to the right of the $k$-th
column where $k=\max(s,s')$.

Since $B^{r,s}\in\CC(A_{N-1}^{(1)})$ is irreducible as a
$U_q(A_{N-1})$-crystal, we have:
\begin{equation} \label{eq:AD}
  \Db_{B^{r,s}}=0 \qquad\text{in type $A_{N-1}^{(1)}$.}
\end{equation}

\subsection{Dual KR crystals of type $A$}
It follows from \cite{S02} that there is a
$U'_q(A_{N-1}^{(1)})$-crystal isomorphism
\begin{equation} \label{eq:AKRdual}
B^{r,s\vee}\cong B^{N-r,s}.
\end{equation}
Let $B^1=B^{1,1}\in\CC(A_{N-1}^{(1)})$. Its crystal graph is
pictured in Figure \ref{fig:crystals} except that $n$ should be
replaced by $N-1$. As ordered sets with respect to the partial order
given in section \ref{sss:boxKR} we have
\begin{align}
B^1&=\{1<2<\dotsm<N\} \\
\label{eq:Bvorder} B^{1\vee}&=\{N^\vee <\dotsm < 2^\vee < 1^\vee \}
\end{align}
As in section \ref{sss:rowKR} $B^{s\vee}$ is the set of weakly
increasing words of length $s$ with letters in $B^{1\vee}$, and its
$U_q(A_{N-1})$-crystal structure is determined from that of
$B^{1\vee}$ using a version of the embedding $\iota$ of
\eqref{eq:rowtensor}. As a special case of \eqref{eq:AD},
\begin{equation} \label{eq:Ddual}
  \Db_{B^{s\vee}}=0
\end{equation}

\begin{prop} \label{pp:Rvee} The category $\CC=\CC(A_{N-1}^{(1)})$
is closed under $\vee$. Moreover, if $B_1,B_2\in
\CC=\CC(A_{N-1}^{(1)})$, then the following diagram commutes:
\begin{equation}
\begin{CD}
  B_2 \otimes B_1 @>{R}>> B_1 \otimes B_2 \\
  @V{\vee}VV @VV{\vee}V \\
  B_1^\vee \otimes B_2^\vee @>>{R}> B_2^\vee \otimes B_1^\vee
\end{CD}
\end{equation}
where $\vee:B_2\otimes B_1\rightarrow B_1^\vee \otimes B_2^\vee$
sends $b_2\otimes b_1$ to $b_1^\vee \otimes b_2^\vee$.
\end{prop}
\begin{proof} $\CC$ is closed under duals by \eqref{eq:AKRdual}.
The composition $b_1^\vee\otimes b_2^\vee \mapsto b_2\otimes b_1
\mapsto R(b_2\otimes b_1):=b_1'\otimes b_2' \mapsto b_2^{\prime\vee}
\otimes b_1^{\prime\vee}$ is a $U'_q(A_{N-1}^{(1)})$-crystal
isomorphism $B_1^\vee \otimes B_2^\vee\rightarrow B_2^\vee \otimes
B_1^\vee$ by \eqref{eq:veeraise}. The commutation follows by the
uniqueness of the $R$-matrix.
\end{proof}

Using section \ref{ss:KRA} we calculate some combinatorial
$R$-matrices and local coenergy functions explicitly.

\begin{prop} \label{pp:RBBv}
\begin{align}
\label{eq:RBBv}
  R_{B^1,B^{1\vee}}(i \otimes j^\vee) &=
  \begin{cases}
   j^\vee \otimes i & \text{if $i\not=j$} \\
    (i+1)^\vee \otimes (i+1) & \text{if $i=j<N$} \\
    1^\vee \otimes 1 & \text{if $i=j=N$.}
  \end{cases} \\
\label{eq:RBvB}
  R_{B^{1\vee},B^1}(j^\vee \otimes i) &=
  \begin{cases}
   i \otimes j^\vee & \text{if $i\not=j$} \\
    (i-1)^\vee \otimes (i-1)^\vee & \text{if $i=j>1$} \\
    N\otimes N^\vee & \text{if $i=j=1$.}
  \end{cases} \\
\label{eq:HBBv}
  H_{B^1 \otimes B^{1\vee}}(i \otimes j^\vee) &=
  \begin{cases}
  1 & \text{if $i=j=2n$} \\
  0 &\text{otherwise}
  \end{cases} \\
\label{eq:HBvB}
  H_{B^{1\vee} \otimes B^1}(j^\vee \otimes i) &=
  \begin{cases}
  1 & \text{if $i=j=1$} \\
  0 &\text{otherwise.}
  \end{cases}
\end{align}
\end{prop}

By Proposition \ref{pp:Hbraided}, the $R$-matrix
$R_{q,p^\vee}:=R_{B^{1\otimes q},B^{1\vee\otimes p}}$ may be
computed using $R_{1,1^\vee}=R_{B^1,B^{1\vee}}$ repeatedly.

\begin{prop} \label{pp:rowRBBv}
Let $b_2\otimes b_1\in B^t \otimes B^{s\vee}\in\CC(A_{N-1}^{(1)})$,
$m=\min(m_N(b_2),m_{N^\vee}(b_1))$, $b_2=\hat{b_2}N^m$, $b_1 =
N^{\vee m} \hat{b_1}$, and $\hat{b_1}'\otimes \hat{b_2}' =
R_{t-m,(s-m)^\vee}(\hat{b_2}\otimes \hat{b_1})$. Then $R_{B^t
\otimes B^{s\vee}}(b_2\otimes b_1) = \hat{b_1}'1^{\vee m} \otimes
1^m \hat{b_2}'$.
\end{prop}

\begin{prop} \label{pp:Rhalf} Let $b_2\otimes b_1\in B^t \otimes
B^{s\vee}$ and $b_1'\otimes b_2' = R_{B^t,B^{s\vee}}(b_2\otimes
b_1)$. Then each $i^\vee \in b_1$ either is unchanged in passing to
$b_1'$ or it changes into $(i+1)^\vee$, including the case of
$N^\vee$ staying the same or changing to $1^\vee$. Similarly each
$i\in b_2$ either remains an $i$ in $b_2'$ or changes to an $(i+1)$
in $b_2'$.
\end{prop}
\begin{proof} This follows from the description of
$R_{B^t,B^{s\vee}}$, \eqref{eq:RBBv}, and the fact that $B^t$ and
$B^{s\vee}$ consist of weakly increasing words.
\end{proof}

\begin{prop} \label{pp:Knuthdual} Let $b,c\in B^{1\vee\otimes m}$
with the total ordering on $B^{1\vee}$ given by \eqref{eq:Bvorder}.
Then $b\Knuth c$ if and only if $b^\vee \Knuth c^\vee$.
\end{prop}

\subsection{Dynkin reversal} \label{ss:Dynkinrev}
Let $N=2n$, $\gf=A_{N-1}^{(1)}$ and $\gb=A_{N-1}$. There is an
involution on the $U_q(A_{N-1})$-crystal $B(\omega_1)$ given by
$j\mapsto j^*:=N+1-j$. It satisfies
\begin{equation} \label{eq:fstar}
  f_i(b^*) = e_{N-i}(b)^*\qquad\text{for all $1\le i\le N-1$.}
\end{equation}
Given any $U_q(A_{N-1})$-crystals $B_1,\dotsc,B_m$ which all have
involutions $*:B_i\rightarrow B_i$ satisfying \eqref{eq:fstar}, one
may define the map
\begin{equation} \label{eq:startensor}
\begin{split}
B_m\otimes\dotsm\otimes B_1 &\overset{*}{\rightarrow}
B_1\otimes\dotsm\otimes B_m \\
(b_m\otimes\dotsm \otimes b_1)^* &= b_1^*\otimes\dotsm\otimes b_m^*
\end{split}
\end{equation}
which also satisfies \eqref{eq:fstar}. Taking $B_i=B(\omega_1)$ for
all $i$, one obtains an involution $*$ on $B(\omega_1)^{\otimes m}$.
This restricts to an involution $*$ on $B(m\omega_1)$ (and hence on
$B^m$) via the embedding $\iota_m$ in \eqref{eq:rowtensor}.

Similarly one may define $*$ on $B^{1\vee}$ by
$j^{\vee*}=(N+1-j)^\vee$ and extend the definition of $*$ to
$B^{s\vee}$ as before.

Finally, using \eqref{eq:startensor}, $*$ may be defined on a tensor
product of KR $U_q(A_{N-1}^{(1)})$-crystals of the form $B^t$ and
$B^{s\vee}$; it reverses the order of the tensor factors.

\begin{prop} \label{pp:R*} Let $B_1,B_2\in \CC=\CC(A_{N-1}^{(1)})$.
Then the diagram commutes:
\begin{equation} \label{eq:R*}
\begin{CD}
  B_2 \otimes B_1 @>{R}>> B_1 \otimes B_2 \\
  @V{*}VV @VV{*}V \\
  B_1 \otimes B_2 @>>{R}> B_2\otimes B_1
\end{CD}
\end{equation}
\end{prop}
\begin{proof} The proof is similar to that of Proposition
\ref{pp:Rvee}.
\end{proof}

\section{VXR map}
\label{sec:VXR}

The virtual crystal method of \cite{OSS03} realizes
$U'_q(\gf_n^\se)$-crystals for $\se\in\onetwoset$ by embedding them
into $U'_q(A_{2n-1}^{(1)})$-crystals. The data for the desired map
\eqref{eq:bij} may be obtained using the virtual crystal embedding
followed by an $R$-matrix. The name VXR stands for ``Virtual X and
R-matrix".

\subsection{Virtual crystals for $B^s_\se$} \label{ss:Vrow}
The KR crystals $B_\se^s\in\CC^1(\gf_n^\se)$ were defined in section
\ref{ss:KR}. To compute the $X$ formula explicitly for
$B\in\CC^1(\gf_n^\se)$ we need the $U'_q(\gf)$-crystal structure on
$B_\se^s$ \cite{KKM94}, the intrinsic coenergy for $B_\se^s$
\cite{HKOTT01,HKOTY99}, and the combinatorial $R$-matrix for
$B_\se^s\otimes B_\se^t$ \cite{HKOT00,HKOT02,NY97}.

Instead of recounting all of these things, we use the method of
virtual crystals \cite{OSS03} to reduce to computations in type $A$.
The virtual X (VX) formula is a formula for the one-dimensional sums
for an arbitrary affine algebra $\gf_X$ expressed in terms of a
simply-laced affine algebra $\gf_Y$. They are known to be valid for
$B\in \CC^1(\gf)$ where $\gf$ is a nonexceptional affine algebra
\cite{OSS03b}. We require the case that $\gf=\gf_n^\se$ for
$\se\in\onetwoset$; see Figure \ref{fig:crystals}. In both cases the
VX formula is stated using the simply-laced affine root system
$A_{2n-1}^{(1)}=A_{N-1}^{(1)}$.

With $\xi$ as in \eqref{eq:xi}, there is an isometric (up to a
scalar multiple) embedding $\Psi$ of the weight lattice of
$\gf_n^\se$ into that of $A_{N-1}^{(1)}$ defined by
\begin{equation} \label{eq:embwt}
\begin{split}
\Psi(\La_i)&=\begin{cases}
    \La_i^A + \La_{N-i}^A & \text{if $i\in\{1,\dotsc,n-1\}$} \\
    \frac{2}{\xi} \La_i^A & \text{if $i\in\{0,n\}$.}
\end{cases} \\
\Psi(\delta) &= \dfrac{2}{\xi} \,\delta^A
\end{split}
\end{equation}
where $\{\La_i\mid 0\le i\le n\}$ and $\delta$ are the fundamental
weights and null root of $\gf_n^\se$ and $\{\La_i^A\mid 0\le i\le
N-1 \}$ and $\delta^A$ are those of $A_{N-1}^{(1)}$. From this it
follows that
\begin{equation}\label{eq:embroot}
\Psi(\alpha_i)=\begin{cases}
    \alpha_i^A + \alpha_{N-i}^A & \text{if $i\in\{1,\dotsc,n-1\}$} \\
    \frac{2}{\xi} \alpha_i^A & \text{if $i\in\{0,n\}$.}
\end{cases}
\end{equation}
There is an induced map between the weight lattices of $\gb_n^\se$
and $A_{N-1}$.

Define the virtual lowering operators to be the composite operators
on a crystal of type $U'_q(A_{N-1}^{(1)})$ defined by
\begin{equation*}
\fh_i = \begin{cases}
  f_i f_{N-i} &\text{if $i\in\{1,2,\dotsc,n-1\}$} \\
   f_i^{\frac{2}{\xi}} & \text{if $i\in\{0,n\}$.}
\end{cases}
\end{equation*}
Define the virtual raising operators $\eh_i$ with $e$ replacing $f$
in the above definition.

For our purposes, an \textbf{aligned virtual crystal} of type $\se$
is a pair $(V,\Vh)$ where $\Vh$ is a tensor product of KR crystals
for $U'_q(A_{2n-1}^{(1)})$, $V\subset \Vh$ is a subset containing
$u(\Vh)$ that is closed under the virtual operators $\fh_i$ and
$\eh_i$ for $i\in \{0,\dotsc,n\}$, such that there is a tensor
product $B$ of KR crystals of type $\se$ and a bijection
$\Psi:B\rightarrow V$ such that
\begin{enumerate}
\item $\Psi(u(B))=\Psi(u(\Vh))$.
\item For all $b\in B$, $\Psi(f_i b)=\fh_i \Psi(b)$
and $\Psi(e_i b)=\eh_i \Psi(b)$.
\end{enumerate}
Moreover, the following must also hold for all $b\in B$, $i\in
\{1,\dotsc,n-1\}$, and $j\in \{0,n\}$:
\begin{equation} \label{eq:aligned}
\begin{split}
\vp_i(b)&=\vp_i(\Psi(b))=\vp_{N-i}(\Psi(b)) \\
\ve_i(b)&=\ve_i(\Psi(b))=\ve_{N-i}(\Psi(b)) \\
\frac{2}{\xi} \, \vp_j(b)&=\vp_j(\Psi(b)) \\
\frac{2}{\xi} \, \ve_j(b)&=\ve_j(\Psi(b))
\end{split}
\end{equation}
From this it follows that
\begin{equation} \label{eq:Psiwt}
\wt(\Psi(b))=\Psi(\wt(b)) \qquad\text{for all $b\in B$.}
\end{equation}
Define $\Vh^s,\Vh^{s\prime}\in \CC(A_{N-1}^{(1)})$ by
\begin{equation} \label{eq:rowV}
\Vh^s = B^{s\vee} \otimes B^s \qquad\qquad \Vh^{s\prime} = B^s
\otimes B^{s\vee}
\end{equation}
and maps $\Psi_s:B_\se^s\rightarrow \Vh^s$ and
$\Psi_s':B_\se^s\rightarrow \Vh^{s\prime}$ by:
\begin{equation} \label{eq:Vrow}
\begin{aligned}
 \Psi_s(u_{B_\se^s}) &= u(\Vh^s) = N^{\vee s} \otimes 1^s &\qquad \Psi_s'(u_{B_\se^s}) &= u(\Vh^{s\prime})
    = 1^s \otimes N^{\vee s}  \\
 \Psi_s(f_i(b)) &= \fh_i \Psi_s(b) &\qquad \Psi_s'(f_i(b)) &= \fh_i \Psi_s'(b)
\end{aligned}
\end{equation}
for all $0\le i\le n$ and $b\in B_\se^s$. By \cite{OSS03b} $\Psi_s$
and $\Psi_s'$ are well-defined injections. Let $V^s=\Image(\Psi_s)$
and $V^{s\prime}=\Image(\Psi_s')$. Then $(V^s,\Vh^s)$ and
$(V^{s\prime},\Vh^{s\prime})$ are aligned virtual crystals with
bijections $\Psi_s:B_\se^s \rightarrow V^s$ and $\Psi_s':B_\se^s
\rightarrow V^{s\prime}$. The map $\Psi_s$ and set $V^s$ are
explicitly given in \cite{OSS03b}. One also has
\begin{equation} \label{eq:PsiR}
  \Psi_s' = R_{B^{s\vee},B^s} \circ \Psi_s.
\end{equation}
On classical highest weight vectors one has
\begin{equation} \label{eq:embhwv}
  \Psi(1^r) = N^{\vee r} 1^{\vee(s-r)} \otimes 1^s
\end{equation}
for $r$ as in \eqref{eq:Bdecomp}.

\begin{example} \label{ex:PsiBox}
For $s=1$ the maps $\Psi$ and $\Psi'$ are given explicitly by
\begin{equation} \label{eq:D2toA}
\begin{aligned}
  \Psi(i) &= i^{*\vee} \otimes i &\qquad \Psi'(i) &= i\otimes i^{*\vee} \\
  \Psi(\bar{i}) &= i^\vee \otimes i^* &\qquad \Psi'(\bar{i}) &= i^*\otimes i^\vee \\
  \Psi(0) &= (n+1)^\vee \otimes (n+1) &\qquad \Psi'(0) &= n \otimes n^\vee \\
  \Psi(\nn) &= 1^\vee \otimes 1 &\qquad \Psi'(\nn) &= N \otimes
  N^\vee
\end{aligned}
\end{equation}
for $1\le i\le n$ where $i^*$ and $i^{*\vee}=i^{\vee*}$ are defined
in section \ref{ss:Dynkinrev}.
\end{example}

Aligned virtual crystals form a tensor category \cite{OSS03}: if
$(V_1,\Vh_1)$ and $(V_2,\Vh_2)$ are aligned virtual crystals of type
$\se$ with bijections $\Psi_i:B_i\rightarrow V_i$, then $(V_2\otimes
V_1,\Vh_2\otimes \Vh_1)$ is an aligned virtual crystal of type $\se$
with bijection $\Psi_2\otimes \Psi_1:B_2\otimes B_1 \rightarrow V_2
\otimes V_1$.

The virtual crystal for an arbitrary $B\in\CC$ is obtained by
tensoring together the virtual crystals for its tensor factors. We
introduce notation for such virtual crystals. Let
$\nu=(\nu_1,\nu_2,\dotsc,\nu_m)$ be a sequence of positive integers.
Define $B_\se^\nu\in \CC^1(\gf_n^\se)$ and
$\Vh^\nu,\Vh^{\nu\prime}\in\CC(A_{2n-1}^{(1)})$ by
\begin{equation} \label{eq:VB}
\begin{split}
B_\se^\nu&=B_\se^{\nu_m}\otimes\dotsm\otimes B_\se^{\nu_1}  \\
\Vh^\nu &= \Vh^{\nu_m} \otimes\dotsm\otimes \Vh^{\nu_1}
\\
\Vh^{\nu\prime} &= \Vh^{\nu_m\prime} \otimes\dotsm\otimes
\Vh^{\nu_1\prime}
\end{split}
\end{equation}
Then $\Psi_\nu:B_\se^\nu\rightarrow \Vh^\nu$ defined by
$\Psi_\nu=\Psi_{\nu_m} \otimes\dotsm\otimes \Psi_{\nu_1}$ gives the
virtual crystal $(V^\nu,\Vh^\nu)$ where $V^\nu=\Image(\Psi_\nu)$.
$\Psi_\nu':B_\se^\nu\rightarrow \Vh^{\nu\prime}$ is defined
similarly. Again
\begin{equation} \label{eq:PsinuR}
  \Psi_\nu' = R \circ \Psi_\nu
\end{equation}
where $R$ is the tensor product of $R_{B^{\mu_i\vee},B^{\mu_i}}$.

Let $L$ be the multiplicity array for $B_\se^\nu$, that is,
$L_i^{(1)}$ is the number of times the part $i$ occurs in $\nu$ and
$L_i^{(a)}=0$ for $a>1$. For $\la\in\Pb^+$, define
$\Pv(B_\se^\nu,\la)$ to be the set of classical highest weight
vectors in $V^\nu\subset \Vh^\nu$ of weight $\Psi(\la)$. The VX
formula is given by
\begin{equation*}
\VX(B_\se^\nu,\la)(t) = \sum_{b\in \Pv(B_\se^\nu,\la)}
t^{\frac{1}{\xi}\Db_{\Vh^\nu}(b)}.
\end{equation*}
The coenergy function is computed in the
$U'_q(A_{N-1}^{(1)})$-crystal $\Vh^\nu$.

\begin{theorem} \cite{OSS03b} For $B_\se^\nu\in \CC^1(\gf_n^\se)$
and $L$ as above and for $\la\in\Pb^+$,
$$X_{L,\la}^\se(t) = \VX(B_\se^\nu,\la)(t).$$
\end{theorem}

For later we need the following fact: the image of $\Psi_\nu$
consists of elements that are fixed by the map $d\mapsto d^{\vee*}$
up to $R$-matrices.

\begin{prop} \label{pp:Psiselfdual} \cite{OSS03} For all
$b\in B_\se^\nu$,
\begin{equation} \label{eq:Psiselfdual}
  \Psi_\nu'(b) = \Psi_\nu(b)^{\vee*}.
\end{equation}
\end{prop}
\begin{proof} One may reduce to the case of a single tensor factor using
\eqref{eq:veetensor} and \eqref{eq:startensor}. For a single tensor
factor \eqref{eq:Psiselfdual} follows from \eqref{eq:RBBv}.
\end{proof}

The $R$-matrix in type $\se\in\onetwoset$ is achieved by a
composition of $R$-matrices of type $\nn$ in the virtual crystal.
Suppose $(V_1,\Vh_1)$ and $(V_2,\Vh_2)$ are aligned virtual crystals
of type $\se$ with bijections $\Psi_i:B_i\rightarrow V_i$. The
combinatorial $R$-matrix $R_{\Vh_2,\Vh_1}:\Vh_2\otimes
\Vh_1\rightarrow \Vh_1\otimes \Vh_2$, restricts to a bijection
called the \textbf{virtual combinatorial $R$-matrix} $\Rv:V_2\otimes
V_1\rightarrow V_1\otimes V_2$ such that the diagram commutes
\cite{ScS04}:
\begin{equation} \label{eq:VR}
\begin{CD}
B_2\otimes B_1 @>R>> B_1 \otimes B_2 \\
@V{\Psi_2\otimes \Psi_1}VV @VV{\Psi_1\otimes \Psi_2}V \\
V_2 \otimes V_1 @>>{\Rv}> V_1 \otimes V_2
\end{CD}
\end{equation}
This is a bijective realization of the invariance of the $\VX$
formula under permutations of the parts of $\nu$.

\subsection{The VXR map}
For a composition $\nu$ we define a map $\eta_\nu$ as in
\eqref{eq:bij}. Let $L$ be the multiplicity array for $B_\se^\nu\in
\CC^1(\gf_n^\se)$. Define $B^\nu,B_\vee^\nu \in \CC(A_{N-1}^{(1)})$
by
\begin{equation} \label{eq:BBvee}
\begin{split}
B^\nu&=B^{\nu_m}\otimes\dotsm\otimes B^{\nu_1} \\
B_\vee^\nu &= B^{\nu_m \vee} \otimes\dotsm\otimes B^{\nu_1\vee}
\end{split}
\end{equation}
We have a commutative diagram
\begin{equation} \label{eq:VXRrow}
\xymatrix{
  {B_\se^\nu} \ar[d]_{\Psi_\nu} \ar[r]^= &
    {B_\se^\nu} \ar[d]^{\Psi_\nu'} \\
  {\Vh^\nu} \ar[d]_{R_+} \ar[r]^{R_0} &
 {\Vh^{\nu\prime}} \ar[d]^{R_-} \\
  {B_\vee^\nu \otimes B^\nu} \ar[r]_R &
  {B^\nu \otimes B_\vee^\nu}
}
\end{equation}
where $R_+^{-1},R_0,R_-,R$ are all compositions of $R$-matrices of
the form $R_{ B^{s\vee} \otimes B^t }$ applied at adjacent tensor
positions. Let $b\in P_\se(B_\se^\nu,\la)$ and
\begin{equation} \label{eq:dcdef}
\dck \otimes c = R_+(\Psi_\nu(b)) \qquad d \otimes \cck =
R_-(\Psi_\nu'(b))
\end{equation}

\begin{example} \label{ex:dckc} In the running example it is convenient to take
$n=3$ even though this is smaller than the generally necessary bound
for $n$. With $b$ as in Example \ref{ex:X},
\begin{equation*}
\begin{tabular}{ccccccccc}
$b$ & = & $\bb$ & 1 & $\nn$ & $\ba$ & 2 & 1 & 1 \\
$\Psi(b)$ & = & $\cb5$ & $\cf1$ & $\ca1$ & $\ca6$ & $\ce2$ & $\cf1$
& $\cf1$ \\
$\Psi'(b)$ & = & $5\cb$ & $1\cf$ & $6\cf$ & $6\ca$ & $2\ce$ & $1\cf$
& $1\cf$
\end{tabular}
\end{equation*}
The map $R_+$ sending $\Psi(b)$ to $\dck \otimes c$ is computed
below. The first line is $\Psi(b)$. To obtain each successive line,
the $R$-matrix \eqref{eq:RBBv} is applied to each pair $i \otimes
j^\vee$. The last line is $\dck\otimes c$.
\begin{equation*}
\begin{array}{cccccccccccccc}
\twocol{7}&\twocol{6}& \twocol{5} & \twocol{4} & \twocol{3} & \twocol{2} & \twocol{1} \\ \cline{1-14} \\[-8pt]
\cb&5&\cf&1&\ul{\ca}&\ul{1}&\ca&6&\ce&2&\cf&1&\cf&1 \\
\cb&\cf&5&\cb&2&\cb&2&\ce&\ul{6}&\ul{\cf}&2&\cf&1&1 \\
\cb&\cf&\cb&5&\cc&3&\ce&2&\ul{\ca}&\ul{1}&\cf&2&1&1 \\
\cb&\cf&\cb&\cc&5&\ce&3&\ca&2&\cf&1&2&1&1 \\
\cb&\cf&\cb&\cc&\cf&6&\ca&3&\cf&2&1&2&1&1 \\
\cb&\cf&\cb&\cc&\cf&\ca&\ul{6}&\ul{\cf}&3&2&1&2&1&1 \\
\cb&\cf&\cb&\cc&\cf&\ca&\ul{\ca}&\ul{1}&3&2&1&2&1&1
\end{array}
\end{equation*}
Therefore $\dck=\cb\cf\cb\cc\cf\ca\ca$ and $c=1321211$.
\end{example}

\begin{lemma} \label{lem:classres} $\dck \otimes c$ and $c$ are
classical highest weight vectors.
\end{lemma}
\begin{proof} Since $b$ is a highest weight vector, $\ve_i(b)=0$ for
all $1\le i\le n$. By \eqref{eq:aligned} $\ve_i(\Psi(b))=0$ for
$1\le i\le N-1$, that is, $\Psi(b)$ is a classical highest weight
vector. Since $R_+$ is an affine crystal isomorphism, $\dck\otimes c
= R_+(\Psi(b))$ is also a classical highest weight vector. By
\eqref{eq:ftensor}, $c$ is a classical highest weight vector.
\end{proof}
Let
\begin{equation} \label{eq:Rtaudef}
  \tau=\wt(c)
\end{equation}
\begin{equation} \label{eq:Bhalf}
\begin{split}
  B^{1\vee}_- &=\{N^\vee<\dotsm<(n+1)^\vee\} \\
  B^{1\vee}_+ &=\{n^\vee<\dotsm<1^\vee \}
\end{split}
\end{equation}
and let $\dcp$ (resp. $\dcm$) be the word obtained by concatenating
the tensor factors in $\dck$ and then restricting to the subalphabet
$B^{1\vee}_+$ (resp. $B^{1\vee}_-$).

Define
\begin{equation} \label{eq:Zmudef}
  Z = [\dcp^\vee] \qquad\qquad \mu=\shape(Z)
\end{equation}

\begin{example} \label{ex:Zmu} For $n=3$, $B^{1\vee}_-=\{\cf<\ce<\cd\}$ and
$B^{1\vee}_+=\{\cc<\cb<\ca\}$. With $\dck$ as in Example
\ref{ex:dckc},
\Yboxdim{10pt}%
\begin{equation*}
\dck = \begin{pmatrix} 7654321 \\ \cb\cf\cb\cc\cf\ca\ca
\end{pmatrix} \quad
\dcp=\begin{pmatrix} 75421 \\ \cb\cb\cc\ca\ca
\end{pmatrix} \quad
\dcm=\begin{pmatrix} 63 \\ \cf\cf
\end{pmatrix} \quad
\dcp^\vee = 11322 \qquad Z=\young(1122,3)
\end{equation*}
\end{example}

The VXR map is defined by
\begin{equation} \label{eq:VXR}
\begin{split}
 \eta_\nu: P_\se(B_\se^\nu,\la) &\rightarrow \bigcup_\tau
 \bigcup_{\mu\in \Par^\se} P_\nn(B^\nu,\tau) \times
 \LR(\tau;\la,\mu) \\
 b &\mapsto (c,Z)
\end{split}
\end{equation}
where the Littlewood-Richardson data $Z$ takes the form of a
$\tau/\la$-Yamanouchi tableau of shape $\mu$; see Theorem
\ref{th:LR}. To prove Theorem \ref{th:KX}, by the discussion in
section \ref{ss:sketch} it suffices to establish the following
result.

\begin{theorem} \label{th:VXR}
$\eta_\nu$ is a well-defined bijection \eqref{eq:bij} satisfying
\eqref{eq:bijstat}.
\end{theorem}

\begin{lemma} \label{lem:details}
\begin{enumerate}
\item $\tau\supset \la$.
\item $\dcm$ is a highest weight vector of weight
$-w_0\la=(0^n,-\la_n,\dotsc,-\la_1)$ where $w_0$ is the longest
element in the Weyl group of type $A_{N-1}$.
\item The word $\dcp^\vee$ is $\tau/\la$-Yamanouchi.
\end{enumerate}
In particular $c\in P_\nn(B^\nu,\tau)$ for some partition $\tau$ and
$Z$ is a $\tau/\la$-Yamanouchi tableau of some partition shape
$\mu$.
\end{lemma}
\begin{proof} By the definition of $\Psi$ and the fact that
$R$-matrices are isomorphisms,
\begin{equation} \label{eq:embedwt} \wt(\dck\otimes c) =
\wt(\Psi(b)) =\la-w_0\la = (\la_1,\dotsc,\la_n,-\la_n,\dotsc,-\la_1)
\end{equation}
Since the assertions to be proved only involve the
$U_q(A_{N-1})$-crystal structure, using the $U_q(A_{N_1})$-crystal
embeddings $B^s\rightarrow B^{1\otimes s}$ and $B^{s\vee}\rightarrow
B^{1\vee\otimes s}$ of \eqref{eq:rowtensor} we may regard
$\dck\otimes c$ as a word in the alphabet $B^1\cup B^{1\vee}$.
Consider the sequence of weights of the right factors of
$\dck\otimes c$ as it is scanned from right to left, letter by
letter. They are all dominant weights by Lemma \ref{lem:classres},
starting with the zero weight and ending at $\la-w_0\la$. Consider
the first $n$ parts of the changing weights. They go from $(0^n)$ to
$\tau$, adding $1$ to some part of the weight as $c$ is scanned from
right to left, and then to $\la$ by removing cells corresponding to
letters in $\dcp$. Consider the last $n$ parts of the weights. They
remain fixed at $(0^n)$ as $c$ is scanned, and then go to
$(-\la_n,\dotsc,-\la_1)$, subtracting $1$ from some part of the
weight, for each letter in $\dcm$. This proves the various
assertions.
\end{proof}

\begin{remark} The splitting of all of the factors $B^s$ into $B^1$
that was used in the above proof, preserves classical crystal
structure but destroys the affine structure and therefore the
coenergy. It turns out that splitting the rightmost tensor factor
preserves the coenergy even though it does not respect the affine
crystal structure. This is exploited in the next section.
\end{remark}

Theorem \ref{th:VXR} is proved by descending induction on the number
of parts of $\nu$. We shall reduce to the case $\nu=(1,1,\dotsc,1)$
using a system of crystal embeddings.

\subsection{Right splitting}
\label{ss:rs}

\begin{prop} \label{pp:split} \cite{S02a,ScS04} Let $\CC$ be either
$\CC(A_n^{(1)})$ or $\CC^1(\gf_n^\se)$ and $B^{r,s}\in \CC$. Then
there is a $U_q(\gb)$-crystal embedding
\begin{equation} \label{eq:rs}
 \rs_{r,s}: B^{r,s} \rightarrow B^{r,s-1}\otimes B^{r,1}
\end{equation}
called \textbf{right splitting}, such that, for any $B\in \CC$, the
map $1 \otimes \rs_{r,s}:B \otimes B^{r,s} \rightarrow B \otimes
B^{r,s-1} \otimes B^{r,1}$ is a $U_q(\gb)$-crystal embedding that
preserves intrinsic coenergy.
\end{prop}

We need explicit formulas for $\rs$ \cite{ScS04}.

For $A_{N-1}^{(1)}$, we write $\rs_s:B^s\rightarrow B^{s-1}\otimes
B^1$ for $\rs_{1,s}$, and $\rs_{s\vee}: B^{s\vee}\rightarrow
B^{(s-1)\vee} \otimes B^{1\vee}$ for $\rs_{N-1,s}$; see
\eqref{eq:AKRdual}. Explicitly $\rs_s(b_s\dotsm b_1)=b_s\dotsm b_2
\otimes b_1$ where $b_i\in B^1$. The same rule works for
$\rs_{s\vee}$ except that $b_i\in B^{1\vee}$. These maps literally
split off the rightmost letter.

For $\gf_n^\se$ we write $\rs^\se_s:B^s_\se \rightarrow B^{s-1}_\se
\otimes B^1_\se$ for $\rs_{1,s}$. Explicitly, if $b\in
B(r\omega_1)\subset B^s_\se$ for $r\ge1$, then one may write $b=ux$
with $x\in B(\omega_1)$ and $u\in B((r-1)\omega_1)$. In this case
$\rs_s^\se(b)=u\otimes x$. Otherwise $r=0$ and $b=\nn\in B(0)\subset
B^s_\se$ and $\rs_s^\se(b)=\bar{1}\otimes 1$.

\begin{theorem} \label{th:Vsplit} \cite{ScS04} Let
$B_\se\in\CC^1(\gf_n^\se)$ have the virtual crystal $(V,\Vh)$ with
bijection $\Psi_B:B_\se\rightarrow V$. Then the following diagram
commutes:
\begin{equation} \label{eq:Vsplit}
\begin{CD}
  B_\se \otimes B_\se^s @>{\Psi_B\otimes \Psi_s}>> V \otimes V^s
  \\
  @V{1\otimes \rs^\se}VV @VV{1\otimes \rs'}V \\
  B_\se \otimes B_\se^{s-1} \otimes B_\se^1 @>>{\Psi_B\otimes \Psi_{s-1}\otimes \Psi_1}> V \otimes V^{s-1}
  \otimes V^1
\end{CD}
\end{equation}
where $\rs':V^s\rightarrow V^{s-1} \otimes V$ is a composition of
$R$-matrices and embeddings $\rs_\nn$.
\end{theorem}

\begin{example} \label{ex:rs'}
The map $\rs':V^s\rightarrow V^{s-1} \otimes V^1$ is given by the
restriction to $V^s$ of the map $\rsh:\Vh^s\rightarrow \Vh^{s-1}
\otimes \Vh^1$ defined by \cite{ScS04}
\begin{equation*}
\xymatrix{ {B^{s\vee} \otimes B^s} \ar[r]^-{R} &{B^s\otimes
B^{s\vee}} \ar[r]^-{1\otimes \rs_{s\vee}} & {B^s\otimes
B^{(s-1)\vee} \otimes B^{1\vee}} \ar[r]^-{R} &}
\end{equation*}
\begin{equation*}
\xymatrix{ {B^{(s-1)\vee} \otimes B^\vee \otimes B^s}
\ar[r]^-{1\otimes 1\otimes \rs_s} &{B^{(s-1)\vee} \otimes B^{1\vee}
\otimes B^{s-1} \otimes B^1} \ar[r]^-{R} & }
\end{equation*}
\begin{equation*}
\xymatrix{ {B^{(s-1)\vee} \otimes B^{s-1} \otimes B^{1\vee}  \otimes
B^1} }
\end{equation*}
\end{example}

\subsection{Reduction from rows to boxes}
We introduce coenergy-preserving embeddings of tensor products of KR
modules into other such tensor products. In type $A$ ($\se=\nn$),
under the identification of $P(L,\la)$ with semistandard tableaux of
shape $\la$ and weight specified by $L$, one may compose such
embeddings and combinatorial $R$-matrices to define an embedding
which coincides with the map of Lascoux \cite{L91}, which sends a
semistandard tableau to a standard one of the same shape and
cocharge. The generalization of this embedding for a general tensor
product of KR crystals in type $A_n^{(1)}$ is given in \cite{S02a}.
Using virtual crystals this can be generalized to the affine types
that embed in type $A^{(1)}$.

Let $\nu=(\nu_1,\dotsc,\nu_m)$ be a sequence of positive integers.
One may pass from $\nu$ to $(1^{|\nu|})$ using two operations:
\begin{enumerate}
\item Exchanging adjacent parts, say, the $r$-th and $r+1$-th, resulting in
$s_r \nu$.
\item Passing to $\rs(\nu) := (1,\nu_1-1,\nu_2,\dotsc,\nu_m)$ if
$\nu_1 \ge 2$.
\end{enumerate}
Define the map
\begin{equation}
  \thet{\se}{\nu}{s_r\nu}:
  P_\se(B_\se^\nu,\la) \rightarrow P_\se(B_\se^{s_r \nu},\la)
\end{equation}
by an $R$-matrix acting in the $r$-th and $(r+1)$-th positions. This
is a weight- and grade-preserving bijection; see Remark
\ref{rem:Border}. If $\nu_1\ge2$ then define
\begin{equation}
  \thet{\se}{\nu}{\rs(\nu)}: P_\se(B_\se^\nu,\la) \rightarrow
  P_\se(B_\se^{\rs(\nu)},\la)
\end{equation}
by $\rs$ acting on the rightmost tensor factor. This is a weight-
and grade-preserving embedding by Theorem \ref{th:Vsplit}.

We reduce the proof of Theorem \ref{th:VXR} to the case
$\nu=(1,1,\dotsc,1)$. The large rank assumption is used here in an
essential way. See also Remark \ref{rem:tauparts}.

\begin{prop} \label{pp:nu} Let $\rho=s_r\nu$ or $\rho=\rs(\nu)$ if $\nu_1\ge2$.
Suppose $\eta_\rho$ is a bijection. Then
$\eta_\nu$ is a bijection such that the diagram commutes:
\begin{equation} \label{eq:theta}
\begin{CD}
  P_\se(B_\se^\nu,\la) @>{\eta_\nu}>> \displaystyle{\bigcup_\tau
  \bigcup_{\mu\in \Par^\se}
  P_\nn(B^\nu,\tau) \times \LR(\tau;\la,\mu)} \\
  @V{\thet{\se}{\nu}{\rho}}VV @VV{\bigcup\, (\thet{\nn}{\nu}{\rho} \times 1)}V \\
P_\se(B_\se^\rho,\la) @>>{\eta_\rho}> \displaystyle{\bigcup_\tau
\bigcup_{\mu\in \Par^\se} P_\nn(B^\rho,\tau) \times
  \LR(\tau;\la,\mu)}
\end{CD}
\end{equation}
\end{prop}
\begin{proof} The map $\eta_\nu$ is well-defined by
Lemma \ref{lem:details} except for the condition
\begin{equation} \label{eq:seshape}
  \mu\in \Par^\se
\end{equation}
which is nontrivial in the case $\se=\hdom\,$.

Let $b\in P_\se(B_\se^\nu,\la)$, $b'=\thet{\se}{\nu}{\rho}(b)$,
$\dck \otimes c = R_+^\nu(\Psi_\nu(b))$, $\dck' \otimes c' =
R_+^\rho(\Psi_\rho(b'))$, $\eta_\nu(b)=(c,Z)$, and
$\eta_\rho(b')=(c',Z')$. By assumption there is a partition $\tau'$
and a partition $\mu'\in\Par^\se$ such that $c'\in
P_\nn(B^\rho,\tau')$ and $Z'$ is $\tau'/\la$-Yamanouchi of shape
$\mu'$. To show that $\eta_\nu$ is a well-defined injection, it
suffices to prove that
\begin{align}
\label{eq:cc'}
 \thet{\nn}{\nu}{\rho}(c)&=c' \\
\label{eq:ZZ'} Z&=Z'.
\end{align}
Suppose that $\rho=s_r \nu$. Then \eqref{eq:cc'} follows from
\eqref{eq:VR} and the uniqueness of the crystal isomorphisms coming
from compositions of $R$-matrices. Since $\dck$ and $\dck'$ are
related by an $R$-matrix, $\dck\Knuth \dck'$ by Proposition
\ref{pp:KnuthA}. By Proposition \ref{pp:Knuthres} with alphabet
$B^{1\vee}_+$, $\dck|_{B^{1\vee}_+}\Knuth \dck'_{B^{1\vee}_+}$.
Equation \eqref{eq:ZZ'} follows by taking duals and applying
Proposition \ref{pp:Knuthdual} and \eqref{eq:Zmudef}.

Let $\rho=\rs(\nu)$ and $\nu=(s,\nuh)$. The diagram
\begin{equation} \label{eq:diagsplitembed}
\xymatrix{%
{B_\se^{\nuh} \otimes B_\se^s } \ar[r]^{1\otimes \rs^\se}
\ar[d]_{\Psi_\nu} & {B_\se^{\nuh} \otimes B_\se^{s-1} \otimes
B_\se^1} \ar[d]^{\Psi_\rho} \\ %
{\Vh^\nuh\otimes B^{s\vee}\otimes B^s} \ar[r]^(.4){1 \otimes \rsh}
\ar[d]_{R_+^\nuh\otimes 1\otimes1} & {\Vh^\nuh\otimes
B^{(s-1)\vee}\otimes B^{s-1} \otimes B^{1\vee} \otimes B^1}
\ar[d]^{R_+^{\nuh}\otimes1\otimes R_{B^{s-1},B^{1\vee}}\otimes 1} \\ %
{B^{{\nuh}\vee} \otimes B^{\nuh} \otimes B^{s\vee} \otimes B^s}
\ar[r]^(.4){1\otimes 1\otimes \rsh'} \ar[d]_{1\otimes R\otimes 1} &
{B^{\nuh\vee} \otimes B^{\nuh} \otimes B^{(s-1)\vee} \otimes
B^{1\vee} \otimes B^{s-1} \otimes B^1} \ar[d]^{1\otimes R'\otimes 1} \\
{B^{{\nuh}\vee} \otimes B^{s\vee}\otimes B^{\nuh}  \otimes B^s}  &
{B^{{\nuh}\vee} \otimes B^{(s-1)\vee} \otimes B^{1\vee} \otimes
B^{\nuh} \otimes B^{s-1} \otimes B^1} }
\end{equation}
commutes, where $\rsh$ is defined in Example \ref{ex:rs'} and
$\rsh': B^{s\vee} \otimes B^s \rightarrow B^{(s-1)\vee} \otimes
B^{1\vee} \otimes B^{s-1} \otimes B^1$ is defined by $\rsh' = (1
\otimes R_{B^{s-1},B^{1\vee}} \otimes 1) \circ \rsh$. The first
square commutes by Theorem \ref{th:Vsplit} and the second by the
fact that $R^\nuh_+$ and $\rsh$ act on disjoint tensor positions.

The composite map down the left (resp. right) column of
\eqref{eq:diagsplitembed} is $R_+^\nu \circ \Psi_\nu$ (resp.
$R_+^\rho \circ \Psi_\rho$). Let $b\in P_\se(B_\se^\nu)$. One may
write
\begin{equation}
\begin{split}
\dckh \otimes \cht \otimes w \otimes 1^s &\in B^{\nuh\vee} \otimes B^\nuh \otimes B^{s\vee} \otimes B^s \\
\dckh \otimes \cht \otimes v \otimes x \otimes 1^{s-1} \otimes 1
&\in B^{\nuh\vee} \otimes B^\nuh \otimes B^{(s-1)\vee} \otimes
B^{1\vee} \otimes B^{s-1} \otimes B^1
\end{split}
\end{equation}
for the image of $b$ in the sets in the third row of
\eqref{eq:diagsplitembed}. The rightmost factor of $b$ is a
classical highest weight vector in $1^r\in B_\se^s$ with $r$ as in
\eqref{eq:Bdecomp}. By explicit calculation
\begin{align} \label{eq:hwvform}
&\begin{aligned} w &= 1^{\vee s} \\ v \otimes x &= 1^{\vee(s-1)}
\otimes 1^\vee
\end{aligned}
&\qquad&\text{if $r=0$} \\
&\begin{aligned} w &= N^{\vee r} 1^{\vee(s-r)} \\ v \otimes x &=
N^{\vee (r-1)} 1^{\vee(s-r)}\otimes N^\vee
\end{aligned}
&\qquad&\text{if $r\ge1$.}
\end{align}

One may write the images of $b$ in the fourth row of
\eqref{eq:diagsplitembed} in the form
\begin{equation}
\begin{split}
\dckh \otimes w'\otimes \cht' \otimes 1^s &\in B^{\nuh\vee} \otimes B^{s\vee}\otimes B^\nuh  \otimes B^s \\
\dckh \otimes v'\otimes x' \otimes \cht'' \otimes 1^{s-1} \otimes 1
&\in B^{\nuh\vee}  \otimes B^{(s-1)\vee} \otimes B^{1\vee} \otimes
B^\nuh \otimes B^{s-1} \otimes B^1
\end{split}
\end{equation}
To prove \eqref{eq:cc'} and \eqref{eq:ZZ'} it suffices to show that
\begin{equation} \label{eq:thestuff}
\begin{split}
  \cht'&=\cht'' \\
  w'|_{B^{1\vee}_-}&\Knuth v'x'|_{B^{1\vee}_-} \\
  w'|_{B^{1\vee}_+}&\Knuth v'x'|_{B^{1\vee}_+}
\end{split}
\end{equation}
During the $R$-matrices going from the third to fourth rows of
\eqref{eq:diagsplitembed}, in passing from $v\otimes x$ to $v'
\otimes x'$ and from $w$ to $w'$, letters in the subalphabet
$B^{1\vee}_+$ remain in $B^{1\vee}_+$. Some $N^\vee$'s may become
$1^\vee$'s and thereafter remain in $B^{1\vee}_+$; all other letters
in $B^{1\vee}_-$ will remain in $B^{1\vee}_-$. The reason for this
is Proposition \ref{pp:Rhalf} and the large rank assumption; each
letter in $w$ or $v\otimes x$ can only change by one during each
$R$-matrix of the form $R_{B^t,B^{s\vee}}$, and we may assume that
the number $\ell(\nuh)=\ell(\nu)-1$ of tensor factors in $\cht\in
B^{\nuh}$ is strictly less than $n$. Equations \eqref{eq:thestuff}
follow from these considerations and the explicit form of the
highest weight vectors \eqref{eq:hwvform}.

This proves the well-definedness and injectivity of $\eta_\nu$. Note
that the maximum number of actual changes that can occur is at most
the number of nonzero parts of $\tau$; see Remark
\ref{rem:tauparts}.

For the surjectivity of $\eta_\nu$, let $\tau$ be a partition,
$\mu\in \Par^\se$, $c\in P_\nn(B^\nu,\tau)$, and $Z$ a
$\tau/\la$-Yamanouchi tableau of shape $\mu$. Let
$c'=\thet{\nn}{\nu}{\rho}(c)$. By the surjectivity of $\eta_\rho$
there is a $b'\in P_\se(B_\se^\rho,\la)$ such that
$\eta_\rho(b')=(c',Z)$, with $R_+^\rho(\Psi_\rho(b'))=\dck'\otimes
c'$.

Suppose first that $\rho=s_r\nu$. Let $b\in B_\se^\nu$ be such that
$R_j^\se(b)=b'$ where $R_j^\se$ exchanges the $r$-th and $(r+1)$-th
factors in $b$. Let $\eta_\nu(b)=(c'',Z'')$ with
$R_+^\nu(\Psi_\nu(b))=\dck''\otimes c''$. By \eqref{eq:VR}
$R_j^\nn(c'')=c'$ where $R_j^\nn$ exchanges the $r$-th and
$(r+1)$-th tensor factors in $c''$. By definition for $\rho=s_r\nu$,
$\thet{\nn}{\nu}{\rho}=R_j^\nn$, so that $R_j^\nn(c)=c'$. Since
$R_j^\nn$ is an isomorphism $c''=c$. Similarly $\dck''$ and $\dck'$
are related by an $R$-matrix. The desired equality $Z''=Z'$ follows
by Propositions \ref{pp:KnuthA}, \ref{pp:Knuthres}, and
\ref{pp:Knuthdual}.

Otherwise suppose $\nu_1\ge2$ and $\rho=\rs(\nu)$. Let
$\nu=(s,\nuh)$ with $s\ge2$ and $b'=y\otimes v\otimes x\in
B_\se^\rho=B_\se^{\nuh}\otimes B_\se^{s-1}\otimes B_\se^1$. Since
$b'$ is a classical highest weight vector, so is $v\otimes x$. We
list the possibilities for $v\otimes x$, together with the resulting
two rightmost tensor factors in $c'$.
\begin{equation*}
\begin{array}{|c||c|} \hline
v\otimes x & c'_2 \otimes c'_1 \\ \hline \hline %
1^r \otimes 1 & 1^r \otimes 1 \\ \hline %
1^r 2\otimes 1 & 1^r 2\otimes 1 \\ \hline %
\raisebox{-1pt}{\vphantom{L} $1^r \bar{1} \otimes 1$} & 1^r 2\otimes 1 \\ \hline %
1^r \otimes \nn & 1^{r-1}2\otimes 1 \\ \hline %
\nn\otimes \nn & 2\otimes 1 \\ \hline %
\end{array}
\end{equation*}
Now $c'$ is in the image of $\thet{\nn}{\nu}{\rho}=1_{B^{\nuh}}
\otimes \rs_s$. By the definition of $\rs$ in type $A$, it follows
that the only possible $c'$ is $1^r \otimes 1$. Define $b=y \otimes
1^{r+1}\in B_\se^\nu$. It satisfies $\thet{\se}{\nu}{\rho}(b)=b'$.
By the commutativity of the diagram \eqref{eq:theta} and the
injectivity of the vertical maps, it follows that
$\eta_\nu(b)=(c,Z)$ as desired.
\end{proof}

\begin{remark} \label{rem:tauparts}
The proof only requires that $\tau$ have at most $n$ nonzero parts.
\end{remark}

\subsection{Grading}
\label{ss:grading}

We show that the VXR map $\eta_\nu$ satisfies \eqref{eq:bijstat}.
The proof uses the large rank assumption in that it is assumed that
$\tau$ has at most $n$ nonzero parts.

\begin{prop} \label{pp:stat} For
$b\in P_\se(B_\se^\nu,\la)$ and $c$ derived from $b$ by
\eqref{eq:dcdef} the equation \eqref{eq:bijstat} holds.
\end{prop}

\begin{example} \label{ex:stat} Let $L=7$, $\se=\cell\,$, $\tau=(4,2,1)$, and
$\la=(2)$. Let $b\in P_\se(1^L,\la)$ and $c\in P_\nn(1^L,\tau)$ be
as in Example \ref{ex:dckc}. Using Example \ref{ex:X},
\eqref{eq:bijstat} says $2(25/2)=25 = (7 - 2) + 2\cdot 10$.
\end{example}
We have \cite{OSS03}
\begin{equation} \label{eq:Dembed}
  \Db_\nn(\Psi(b)) = 2 \Db_\se(b).
\end{equation}
For the following calculation we assume that $\se=\cell\,$ since
$\se=\hdom$ is essentially a special case.

Since the embeddings $\thet{\se}{\nu}{\rho}$ preserve highest
weights and coenergy, by Proposition \ref{pp:nu} we may reduce to
the case that $\nu=(1,1,\dotsc,1)=1^L$.

Let $B_2=B^{1\vee \otimes L}$ and $B_1=B^{1\otimes L}$. By
Proposition \ref{pp:DR} and \eqref{eq:D2} we have
\begin{equation} \label{eq:Dwork}
  \Db_\nn(\Psi(b)) = \Db_\nn(\dck \otimes c) =
  \Hb_{B_2 \otimes B_1}(\dck \otimes c) +
  \Db_{B_1}(c) + \Db_{B_2}(\cck)
\end{equation}
where $\cck$ is defined by \eqref{eq:dcdef}. To establish
\eqref{eq:bijstat} it suffices by \eqref{eq:Dembed} and
\eqref{eq:Dwork} to show that
\begin{align}
\label{eq:HH} \Hb_{B_2\otimes B_1}(\dck\otimes c) &= L -|\la| \\
\label{eq:Dccheck} \Db_{B_1}(c) &= \Db_{B_2}(\cck).
\end{align}
To verify \eqref{eq:Dccheck} we compute the $R$-matrix $B_2\otimes
B_1\rightarrow B_1 \otimes B_2$ using the composition $R_- \circ R_0
\circ R_+$ where $R_0=R_{B^{1\vee},B^1}^{\otimes L}:(B^{1\vee}
\otimes B^1)^{\otimes L}\rightarrow (B^1 \otimes B^{1\vee})^{\otimes
L}$.

By Propositions \ref{pp:Rvee}, \ref{pp:R*}, and \ref{pp:Psiselfdual}
and equations \eqref{eq:veetensor} and \eqref{eq:startensor}, we
have
\begin{equation*}
  d \otimes \cck = (\dck \otimes c)^{\vee*} = (\dck)^{\vee*}
  \otimes c^{\vee*}.
\end{equation*}
In particular $\cck = c^{\vee*}$. To prove \eqref{eq:Dccheck} it
suffices to prove the following result.

\begin{prop} \label{pp:DAdual} For any $c\in B^{1\otimes L}$,
\begin{equation} \label{eq:DAdual}
  \Db(c) = \Db(c^{\vee*}).
\end{equation}
\end{prop}
\begin{proof} Let $c= c_L\dotsm c_1$ with $c_i\in B^1$. Then $c^{\vee*} =
c_L^{\vee*} \dotsm c_1^{\vee*}$ by \eqref{eq:veetensor} and
\eqref{eq:startensor}. Note that for $i,j\in B^1$, $i<j$ if and only
if $i^{\vee*}<j^{\vee*}$ in $B^{1\vee}$. Computing $\Db(c)$ and
$\Db(c^{\vee*})$ using \eqref{eq:Dn=}, \eqref{eq:DB11},
\eqref{eq:HB11}, and \eqref{eq:Ddual}, equation \eqref{eq:DAdual}
follows.
\end{proof}

To read the rest of this proof it is useful to look at Example
\ref{ex:dckc}. We now prove \eqref{eq:HH}. From the definitions it
is easy to see that
\begin{equation} \label{eq:DLla}
  L-|\la| = 2m + m_\nn(b)
\end{equation}
where $m$ is the number of barred letters in $b$ and $m_\nn(b)$ is
the number of symbols $\nn$ in $b$. By Proposition \ref{pp:Hbraided}
and \eqref{eq:HBBv}, $\Hb_{B_2 \otimes B_1}(\dck\otimes c)$ is equal
to the number of times that the transformation
\begin{equation} \label{eq:trans}
1^\vee \otimes 1 \rightarrow N \otimes N^\vee
\end{equation}
occurs during the application of $R_{B_2\otimes B_1}$ that sends
$\dck\otimes c$ to $d\otimes \cck$. Therefore it suffices to show
that
\begin{enumerate}
\item During $\dck \otimes c\rightarrow \Psi(b)$ \eqref{eq:trans}
occurs $m$ times.
\item During $\Psi(b)\rightarrow \Psi'(b)$ \eqref{eq:trans} occurs
$m_\nn(b)$ times.
\item During $\Psi'(b)\rightarrow d\otimes \cck$ \eqref{eq:trans}
occurs $m$ times.
\end{enumerate}
(2) is immediate. By symmetry (3) follows from (1). For (1) it is
equivalent to show that in passing from $\Psi(b)$ to $\dck\otimes c$
by $R_+^{-1}$ there are $m$ applications of the inverse of
\eqref{eq:trans}. To see this, note that for every barred letter
$b_i$ in $b$ the corresponding plain letter in $\Psi(b_i)$ has value
greater than $n$. As this plain letter is switched to the right by
applications of $R^{-1}$, its value increases until it has value
$N=2n$ and it participates in the inverse of \eqref{eq:trans}. It
must do so, for it ends up as a letter in the highest weight vector
$c$, whose weight $\tau$ is a partition with at most $n$ nonzero
parts. Moreover this letter can only participate once in the inverse
of \eqref{eq:trans} since it switches with at most $n$ dual letters.
Conversely, the plain letter in $\Psi(b_i)$ for $b_i$ unbarred,
being of value at most $n$, cannot increase in value to more than
$N$ when it gets switched to the right, and hence cannot participate
in the inverse of \eqref{eq:trans}.

This proves Proposition \ref{pp:stat}.

\subsection{Incremental VXR} \label{ss:incVXR} For the eventual
proof of the case $\nu=(1,1,\dotsc,1)$ of Theorem \ref{th:VXR} we
describe the change in the various data in passing from $b'$ to $b$
where $b=x \otimes b'$ with $x=b_L\in B_\se^1$ and $b'\in
P_\se(1^{L-1},\la')$. Denote the data corresponding to $b'$ by
primes.

To compare with the DDF bijection, we define tableaux $\Pt,\Qt,\Tt$
as follows. $\Pt$ is the standard tableau associated with $c$. That
is, $c_i=r$ if and only if the letter $i$ is in the $r$-th row of
$\Pt$. Let $\Qt$ be the injective tableau of shape $\tau/\la$ on the
alphabet $A^*$ where $A$ is the set of positions in $\dck$ occupied
by the subword $\dcp$, such that the $i$-th letter of the biword
$\dcp$ is $r^\vee$ if and only if $i^*$ is in the $r$-th row of
$\Qt$. The tableau $\Qt$ is well-defined by Lemma \ref{lem:details}.
Define the tableau $\Tt$ to be the injective tableau of shape $\la$
whose alphabet is the complement of $A$, such that the biletter
$(i^{*\vee},j)$ is in $\dcm$ if and only if $j$ is in the $i$-th row
of $\Tt$. It is well-defined by Lemma \ref{lem:details}.

The map $R_+$ can be computed by the composition
\begin{equation*}
(B^{1\vee} \otimes B^1)^{\otimes L} \rightarrow (B^{1\vee} \otimes
B^1) \otimes B^{1\vee\otimes (L-1)} \otimes B^{1\otimes (L-1)}
\rightarrow (B^{1\vee})^{\otimes L} \otimes B^{1\otimes L}.
\end{equation*}
This sends $\Psi(b) \mapsto \Psi(x) \otimes \dck' \otimes c'\mapsto
\dck \otimes c$. From this it follows that $c=c_L c'$ and hence that
$\Pt'=\Pt|_{<L}$, with extra cell $s$ in the addable cell in the
$c_L$-th row. Let $\Psi(x)=\zv \otimes y$ for some $y,z\in B^1\in
\CC(A_{2n-1}^{(1)})$. We must compute the map $\zv \otimes y \otimes
\dck' \mapsto \dck \otimes c_L$. We see that $\dck$ has last letter
$\zv$. We watch the process that changes $y$ into $c_L$ as it moves
past $\dck'$ by $R$-matrices. By taking sufficiently large $n\ge
\ell(\tau)$ we may assume that $c_L\le n$. Since the moving plain
letter can only interact with $L-1$ dual letters, the moving letter
can change less than $n$ times.

\subsubsection{Suppose $x=r$ for $1\le r\le n$.}
Then $y=r$, $z=r^*=2n+1-r$, and $\la$ is obtained by adding a corner
cell $s'$ to $\la'$ in the $r$-th row. In this case no letters of
$\dcm'$ are changed and $\dcm=z^\vee \dcm'$ where $z^\vee$ has index
$L$. This adjoins the letter $L$ to the $r$-th row of $\Tt'$ to give
$\Tt$. The $R$-matrix computation finds the biletter $(r^\vee,i_r)$
in $\dcp'$ with $i_r$ maximal and changes it into
$((r+1)^\vee,i_r)$. This changes the moving plain letter to $r+1$.
Then it finds the biletter $((r+1)^\vee,i_{r+1})$ in $\dcp'$ with
$i_{r+1}<i_r$ maximal. The algorithm continues in this way, ending
with the plain letter $c_L$ and having changed biletters
$(r^\vee,i_r),\dotsc,((c_L-1)^\vee,i_{c_L-1})$ by changing their
first components from $j^\vee$ to $(j+1)^\vee$. There is no biletter
in $\dcp'$ of the form $((c_L)^\vee,i_{c_L})$ with
$i_{c_L}<i_{c_L-1}$. This is equivalent to the internal insertion
$\Qt=I_{s'}(\Qt')$.

\subsubsection{Suppose $x=\bar{r}$ for $1\le r\le n$.}
Then $z=r$ and $y=r^*=2n+1-r>n$. Since the moving plain letter $y$
is transformed into a letter $c_L\le n$, it follows that there is a
sequence of biletters
$((2n+1-r)^\vee,j_r),((2n-r)^{\vee},j_{r-1}),\dotsc,((2n)^\vee,j_1)$
in $\dcm'$ with $j_r>j_{r-1}>\dotsm>j_1$, each index maximal, and
biletters $(1^\vee,i_1),(2^\vee,i_2),\dotsc$,
$((c_L-1)^\vee,i_{c_L-1})$ in $\dcp'$ with
$j_1>i_1>i_2>\dotsm>i_{c_L-1}$ with indices maximal. These biletters
change such that their first letter $j^\vee$ is sent to
$(j+1)^\vee$, with the exception that
$((2n)^\vee,j_1)\mapsto(1^\vee,j_1)$. This means that $\Tt$ is
obtained from $\Tt'$ by the reverse row insertion at the corner cell
in the $r$-th row, with ejected letter $j_1$. $\dcp$ is obtained by
the above changes in the letters $(j^\vee,i_j)$ in $\dcp'$, with an
extra letter $r^\vee$ at the left end, indexed by $L$. Therefore
$\Qt$ is obtained from $\Qt'$ by the external insertion of $j_1^*$,
together with the adjoining of $L^*$ to the inside at the addable
cell in the $r$-th row of $\la'$.

\subsubsection{Suppose $x=\nn$.} Then $z=y=1$. The moving plain
letter $y=1$ finds the biletter $(1^\vee,j_1)$ with $j_1$ maximal,
then $(2^\vee,j_2)$ with $j_2<j_1$ maximal, etc., ending with
$((c_L-1)^\vee,j_{c_L-1})$. Thus $\Tt=\Tt'$ and $\Qt$ is obtained
from $\Qt'$ by the external insertion of $L^*$.

\section{The bijection via DDF}
\label{sec:DDF}

Let $\se\in \{ \cell\,,\hdom\,\}$ and $\nu=1^L$. This section
defines the DDF map, which is a form of the desired bijection
\eqref{eq:bij} under the present assumptions. The DDF map is the
composition of a number of bijections: Schensted row and column
insertion \cite{Sc61}, the Delest-Dulucq-Favreau (DDF) bijection
\cite{DDF88}, and the Burge correspondence \cite{B74}.

\subsection{Highest weight vectors and sequences of partitions}
We first translate the definition of $P_\se(1^L,\la)$ given in
Example \ref{ex:standardpath} into the language of
multiplicity-space tableaux, which are certain sequences of
partitions.

Recall the encoding of dominant weights by partitions in subsection
\ref{ss:limit}. We shall use this to identify a dominant integral
weight for one of the classical types $A_{n-1}$, $B_n$, $C_n$,
$D_n$, with an element of $\Z^n$.

For each $b\in P_\se(1^L,\la)$ consider the associated sequence of
partitions
$$\lad=(\nn=\la^{(0)}\subset\la^{(1)}\subset\dotsm\subset\la^{(L)}=\la).$$
The sequences of partitions that arise in this manner are
characterized as follows. Every partition $\la^{(i)}$ is obtained
from the previous partition $\la^{(i-1)}$ by
\begin{enumerate}
\item Adding a cell to the $r$-th row if $b_i=r$. This is allowed in all types.
\item Removing a cell from the $r$-th row if $b_i=\bar{r}$.
This is allowed in all types but $\se=\nn$.
\item Doing nothing if $b_i=\nn$. This is allowed only in type $\se=\cell\,$.
\end{enumerate}
These sequences of partitions define three kinds of
multiplicity-space tableaux $\lad$ of shape $\la$ and length $L$.
For type $\nn$ they are called \textbf{standard tableaux}. For types
$\hdom\,,\vdom\,$ they are called \textbf{oscillating tableaux}. For
type $\cell\,$ we call them \textbf{Motzkin tableaux}.

\begin{remark} For large rank, the symbol $0$ does not appear in a
classical highest weight vector of type $\cell\,$. Similarly in type
$\vdom\,$ the symbols $n$ and $\bar{n}$ don't appear. Otherwise
these tableau definitions must be modified.
\end{remark}

\begin{example} \label{ex:hwv} The elements $c\in P_\nn(1^7,(4,2,1))$ and $b\in P_\cell(1^7,(2))$
in Example \ref{ex:X} have corresponding standard and Motzkin
tableaux given by \Yboxdim{5pt}
\begin{align*}
c&=1321211
  &&(\nn, \quad \yng(1), \quad\yng(2),
\quad\yng(2,1),\quad\yng(3,1),\quad\yng(3,2),\quad\yng(3,2,1),\quad\yng(4,2,1))
\\
b&=\bar{2}1\nn\bar{1}211
  &&(\nn,\quad \yng(1), \quad \yng(2),\quad  \yng(2,1),\quad
\yng(1,1),\quad \yng(1,1),\quad  \yng(2,1),\quad\yng(2)).
\end{align*}
\end{example}

\subsection{Robinson-Schensted}
Let $w=w_1w_2\dotsm w_L$ be a word. The row insertion version of the
Robinson-Schensted correspondence \cite{Sc61} is a bijection
$w\mapsto (P,Q)$ between words in an alphabet $A$, and pairs of
tableaux of the same partition shape with $P$ semistandard in the
alphabet $A$ and $Q$ standard. Let $P=[w]$ and let $Q$ be defined by
the sequence of partitions $\lad$ where
$\la^{(i)}=\shape([w_1w_2\dotsm w_i])$. This is denoted
\begin{equation} \label{eq:rowins}
  \bilet{P}{Q} = \left( \bilet{\nn}{\nn} \leftarrow
  \begin{pmatrix}
  1&2&\dotsm&L \\
  w_1&w_2&\dotsm&w_L
  \end{pmatrix} \right).
\end{equation}
The column insertion version of the Robinson-Schensted
correspondence \cite{Sc61} is a bijection between the same sets as
above, sending the word $u=u_L\dotsm u_2 u_1$ to $(P,Q)$ such that
$P=[u]$ and $Q$ is given by the sequence of partitions given by the
shapes of the tableaux $[u_i\dotsm u_2 u_1]$. This is written
\begin{equation} \label{eq:colins}
  \bilet{P}{Q} = \left(
  \begin{pmatrix}
  L&\dotsm&2&1 \\
  u_L&\dotsm&u_2&u_1
  \end{pmatrix} \rightarrow \bilet{\nn}{\nn}
  \right)
\end{equation}
Here the letters of $u$ get inserted started from the \textit{right}
end.

\begin{remark} \label{rem:rowcol} By definition the bumping
(``P'') tableaux for the row and the column insertion of a word $w$,
are equal, as both are defined by $[w]$.
\end{remark}

The Robinson-Schensted correspondence has the following well-known
involution symmetry.

\begin{prop} \label{pp:inv}
Suppose that $w$ is a permutation of the set $\{1,2,\dotsc,L\}$ with
$w(i)$ written $w_i$. Let $w\mapsto (P,Q)$ under the bijection
\eqref{eq:rowins}. Then $w^{-1}\mapsto (Q,P)$ under
\eqref{eq:rowins}. A similar property holds for the bijection
\eqref{eq:colins}.
\end{prop}

\subsection{Skew insertion}
\label{ss:skew} The following skew version of Robinson-Schensted row
insertion is due to Sagan and Stanley \cite{SS90}. Let $T$ be a
tableau of skew shape $\la/\mu$. In this algorithm the pair of
shapes $\mu\subset\la$ is important, not just the set difference of
their partition diagrams.

The \textbf{external} row insertion of a letter $x$ into the skew
tableau $T$ is defined by inserting $x$ into the first row of $T$
and propagating the insertion in the row insertion algorithm. One
may imagine that the shape $\mu$ is occupied by a tableau having all
letters smaller than those in $x$ or $T$. Let $s$ be the newly
created cell. The result is a tableau of skew shape
$(\la\cup\{s\})/\mu$. This process can be reversed given the cell
$s$.

Let $s'$ be an addable corner of $\mu$. The \textbf{internal} row
insertion on $T$ at $s'$ is defined by removing the letter of $T$ at
$s'$ and inserting it into the next row, propagating the insertion
as in the usual row insertion algorithm.

Let $s$ be the newly created cell. If $s'$ is not in $T$ then let
$s=s'$ and declare the resulting tableau to be $T$. In either case
the new pair of shapes is declared to be
$(\la\cup\{s\})/(\mu\cup\{s'\}$. This process can also be reversed
given $s$.

Note that if one is given a tableau $S$ of skew shape $\la/\mu$ and
a corner $s$ for $\la$, the reverse skew row insertion on $S$ at $s$
may be the reverse of either an internal or external insertion,
depending on whether the reverse bumping path ends by hitting the
shape $\mu$ or not.

A skew version of the row insertion Robinson Schensted
correspondence is defined as follows. An \textbf{injective} taleau
is one with no repeated entries. Fix two partitions $\beta$ and
$\gamma$. Let $A$ (values) and $B$ (positions) be alphabets. The
input for the algorithm is a triple $(T,U,w)$ where $w:B'\rightarrow
A'$ is a bijection between subsets $B'\subset B$ and $A'\subset A$,
and $T$ and $U$ are injective tableaux of respective skew shapes
$\beta/\alpha$ and $\gamma/\alpha$ for some $\alpha$, such that the
letters occurring in $T$ and $A'$ (resp. $U$ and $B'$) are disjoint
subsets of $A$ (resp. $B$). The output is a pair $(P,Q)$ of
injective tableaux of respective shapes $\delta/\gamma$ and
$\delta/\beta$ for some $\delta$ with letters taken from the
respective alphabets $A$ and $B$.

At the beginning, set the bumping tableau to be the skew tableau $T$
and the recording tableau to be empty of skew shape $\beta/\beta$.
Let $b$ be the smallest letter among $U$ and $B'$. If $b\in U$ then
perform an internal insertion on the bumping tableau at the cell of
$U$ containing $b$. If $b\in B'$ then externally insert $w(b)$ into
the bumping tableau. Either way, adjoin the letter $b$ to the
recording tableau at the newly created cell. Continue with the
letters of $U$ and $B'$ in increasing order in this manner. Then $P$
is the resulting bumping tableau and $Q$ the recording tableau.

\begin{example} \label{ex:skewins} Let $A=\{1,2,\dotsc,7\}$ and
$B=A^*=\{7^*<6^*<\dotsm<1^*\}$. Let
\begin{equation*}
T=\young(36)\qquad U=\nn \qquad w=\begin{pmatrix} 7^* & 5^* & 4^* & 2^* & 1^*\\
1&5&2&4&7
\end{pmatrix}.
\end{equation*}
Here $A'=\{1,2,4,5,7\}$ and $B'={A'}^*$.
\begin{equation*}
P=\young(1247,35,6) \qquad Q=\young(::\stb\sta,\stg\ste,\std)
\end{equation*}
\end{example}

\subsection{DDF bijection}
Delest, Dulucq, and Favreau \cite{DDF88} defined a bijection from
oscillating tableaux of shape $\la$ and length $L$, to pairs $(T,I)$
where $I$ is a fixed-point-free involution on a subset
$A\subset\{1,2,\dotsc, L\}$ and $T$ is an injective tableau of shape
$\la$ consisting of exactly the letters in the subalphabet
complementary to $A$ in $\{1,2,\dotsc,L\}$. A skew version of the
DDF bijection was given by Dulucq and Sagan \cite{DS95} but we will
not need the extra generality.

We now describe a straightforward extension of the DDF bijection to
a bijection between Motzkin tableaux and pairs $(T,I)$ which are as
above except that the involution $I$ may have fixed points. This
extended DDF bijection also admits a skew version.

Let $b=b_L\dotsm b_2 b_1\in B_\cell(1^L,\la)$ and $\lad$ the
associated Motzkin tableau. Start with $T_0$ the empty tableau and
$I_0$ the empty involution. Suppose $T_{i-1}$ and $I_{i-1}$ have
been constructed.
\begin{enumerate}
\item If $b_i=r$, let $T_i$ be obtained by adjoining the
letter $i$ to $T_{i-1}$ in the $r$-th row and let $I_i=I_{i-1}$.
\item If $b_i=\bar{r}$, let $T_{i-1}$ be the result of the reverse
row insertion on $T_i$ at the corner cell in its $r$-th row. Let $a$
be the ejected letter. Then let $I_i=I_{i-1}(a,i)$, that is, the map
$I_i$ extends $I_{i-1}$ by the 2-cycle $(a, i)$.
\item If $b_i=\nn$, let $T_i=T_{i-1}$ and add the fixed point $(i)$
to the involution $I_{i-1}$ to obtain $I_i$.
\end{enumerate}
The bijection sends $b\rightarrow(T,I)$ where $T=T_L$ and $I=I_L$.

\begin{remark} \label{rem:DDFeven} $b\in B_{\hdom\,}(1^L,\la)$ (that is,
$b$ has no symbols $\nn$) if and only if $I$ has no fixed points
\cite{DDF88}.
\end{remark}

\begin{example} \label{ex:DDF} Let $b$ be as in Example
\ref{ex:hwv}. The computation of the DDF bijection for the Motzkin
tableau associated with $b$ is given in Figure \ref{fig:mDDF} along
with data that will be useful later.
\end{example}
\begin{figure}[ht]
\begin{equation*}
\begin{array}{|c|l|c|l|c|l|l|} \hline
\scalebox{.7}{\vphantom{\young(1,1)}}i&\hphantom{xx}
T&I&\hphantom{xxxx}S&I w_0&\hphantom{xxxx}P&\hphantom{xxxx}Q
\\ \hline\hline
1&\vphantom{\young(1,1)}\young(1) &&  &  & \young(1) & \young(\bl) \\[2mm]
2&\young(12)&&  &  & \young(12) & \young(\bl\bl) \\[2mm]
3&\young(12,3)&&  &  & \young(12,3) & \young(\bl\bl,\bl) \\[4mm]
4&\young(1,3)&(24)& \young(24) & \pmt{\hd\hb}{24} & \young(124,3) & \young(\bl\hd\hb,\bl) \\[4mm]
5&\young(1,3)&(24)(5)& \young(24,5) & \pmt{\he\hd\hb}{524} & \young(124,35) & \young(\bl\he\hb,\bl\hd) \\[4mm]
6&\young(16,3)&(24)(5)& \young(24,5) &\pmt{\he\hd\hb}{524} &  \young(124,35,6) & \young(\bl\bl\hb,\bl\he,\hd) \\[5.5mm]
7&\young(36)&(24)(5)(17)&
\young(1247,5)&\pmt{\hg\he\hd\hb\ha}{15247} &\young(1247,35,6) &
\vphantom{\young(1,2,3,4)} \young(\bl\bl\hb\ha,\hg\he,\hd) \\
\hline\hline
\end{array}
\end{equation*}
\caption{DDF computation} \label{fig:mDDF}
\end{figure}

\subsection{The Burge correspondence}
Column insertion defines a bijection between involutions $I$ on a
set $A$ and standard tableaux $S$ with alphabet $A$. Specifically,
write $I$ as a biword with upper word given by the set $A$ written
in decreasing order, and with $I(a)$ below $a$ for all $a\in A$.
Define the tableau $S$ by
\begin{equation} \label{eq:Sdef}
\bilet{S}{S} = \left( I \rightarrow \bilet{\nn}{\nn} \right),
\end{equation}
which is well-defined by Proposition \ref{pp:inv}. Define
\begin{equation} \label{eq:mudef}
\mu=\shape(S).
\end{equation}

\begin{example} \label{ex:involution}
Let $A=\{1,2,4,5,7\}$ and $I=(2,4)(5)(1,7)$ an involution on $A$.
Then
\begin{equation*}
I= \begin{pmatrix}
  7&5&4&2&1 \\
  1&5&2&4&7
\end{pmatrix} \qquad
S = \young(1247,5)
\end{equation*}
\end{example}

\begin{theorem} \label{th:Burge} \cite{B74} $I$ has no fixed points if and only if
$S$ has even row lengths.
\end{theorem}

\begin{remark} \label{rem:evenBurge} It follows that if $b\in P_\se(1^L,\la)$ for
$\se\in \{\cell\,,\hdom\,\}$ then $\mu\in \Par^\se$.
\end{remark}

Burge gave a direct bijection to compute $S$ from $I$ when $I$ is
fixed-point-free \cite{B74}. Although not necessary for our proof,
we include the following straightforward extension of Burge's
algorithm, which gives the above bijection $I\mapsto S$ for any
involution $I$. First, for each transposition $(a,i)$ in $I$, create
the ordered pair $(a,i)$ where $a<i$. For each fixed point $i$ in
$I$, create the ordered pair $(i,i)$. Sort the ordered pairs
according to second component, giving the sequence
$(a_1,i_1),(a_2,i_2),\dotsc,(a_k,i_k)$, say. Let $S_0$ be the empty
tableau. Given $S_{j-1}$, form $S_j$ as follows.
\begin{enumerate}
\item If $a_j<i_j$ then column insert $a_j$ into $S_{j-1}$,
forming the tableau $S'_j$ with new cell $s$. Let $S_j$ be obtained
from $S'_j$ by adjoining the letter $i_j$ at the bottom of the
column just to the right of that of $s$.
\item If $a_j=i_j$ then $a_j$ is the largest letter so far; adjoin
it at the bottom of the first column of $S_{j-1}$ to obtain $S_j$.
\end{enumerate}

\begin{example} With $I$ as in Example \ref{ex:involution}, we have
the sequence of $(a_j,i_j)$ given by $((2,4),(5,5),(1,7))$. Then
\Yboxdim{10pt}%
\begin{equation*}
S'_1=\young(2)\quad S_1=\young(24)\quad S_2=\young(24,5)\quad
S'_3=\young(124,5)\quad S_3=\young(1247,5)=S
\end{equation*}
\end{example}

\subsection{Insertion of DDF data}
\label{ss:DDFins} Suppose $b\in P_\se(1^L,\la)$ for
$\se\in\{\cell\,,\hdom\,\}$ with DDF data $(T,I)$ and $I\mapsto S$
under the Burge correspondence. Let $A$ be the set on which $I$ is
an involution. For $1\le i\le L$ let $\hat{i}=L+1-i$ and let
$\hat{A}$ be the image of $A\subset \{1,2,\dotsc,L\}$ under
$i\mapsto \hat{i}$. Then let $Iw_0$ denote the bijection from
$\hat{A}$ to $A$ given by $\hat{a}\mapsto I(a)$. Its biword is
obtained from that of $I$ by replacing each upper letter $x$ by
$\hat{x}$.

\begin{example} With $A$ and $I$ as in Example \ref{ex:involution}, the biword for the
bijection $Iw_0$ is given in the last row of Figure \ref{fig:mDDF}.
\end{example}

Define
\begin{equation} \label{eq:DDFPQ}
  \bilet{P}{Q} = \left( \bilet{T}{\nn} \leftarrow (I w_0) \right),
\end{equation}
that is, row insert the biword $Iw_0$ into $T$, and let $P$ and $Q$
be the bumping and recording tableaux. By \eqref{eq:Sdef} and Remark
\ref{rem:rowcol} we have
\begin{equation} \label{eq:DDFP}
  P = [T \cdot S].
\end{equation}
Let $c=c_L\dotsm c_2c_1$ be the Yamanouchi word corresponding to
$P$. Define
\begin{equation}
  \tau=\wt(c).
\end{equation}
Define the DDF map by
\begin{equation} \label{eq:DDF}
\begin{split}
  b\mapsto (c,T\otimes S)
\end{split}
\end{equation}
using Littlewood-Richardson data in the form of Theorem \ref{th:LR}
part (3).

\begin{prop} \label{pp:DDFbij} The DDF map \eqref{eq:DDF} gives a bijection
\eqref{eq:bij}.
\end{prop}
\begin{proof} The map is well-defined by construction and
Remark \ref{rem:evenBurge}, which shows that the tricky condition
$\mu\in \Par^\se$ is satisfied. Since the DDF map was comprised of
bijections, given $(c,T\otimes S)$ one may reconstruct $I$ and then
$b$ by reversing the extended Burge and DDF bijections.
\end{proof}

\begin{example} \label{ex:PTS} Consult the last row of the table in
Figure \ref{fig:mDDF}, which computes $Iw_0$, $P$, and $Q$ from $T$
and $I$ corresponding to $b$ of Example \ref{ex:hwv}. The associated
Yamanouchi word $c$ is given in Example \ref{ex:hwv}. For this data,
$\la=(2)$, $\mu=(4,1)$, and $\tau=(4,2,1)$.
\end{example}

\subsection{Incremental DDF} \label{ss:incDDF}
We describe the change in the output data of the DDF map when one
more letter is added to the path $b$. Let $b=x\otimes b'$ so that
$x=b_L\in B_\se$. Denote by $T',I',S',P',Q'$ the DDF data associated
with $b'$. Let $\la',\mu',\tau'$ be the shapes of $T',S',P'$
respectively.

The biletter $\bilet{a}{i}$ will be written $(a,i)$. Restricting
\eqref{eq:DDFP} to the alphabet of letters less than $L$, by
Proposition \ref{pp:Knuthres}
\begin{equation} \label{eq:Pres}
  P|_{<L} = [(T|_{<L}) (S|_{<L})].
\end{equation}
We will show that in all cases
\begin{equation} \label{eq:Presprime}
  P|_{<L} = P'.
\end{equation}
Then $P$ is obtained by adjoining $L$ to $P'$ at the cell $s$, say.
To determine this cell $s$ it is enough to compute $Q$ from $Q'$.
This is done in each case as follows.
\begin{enumerate}
\item Suppose $x=r$. Let $s'$ be the addable cell in the $r$-th
row of $\la'$. Then $T$ is obtained by adjoining $L$ to $T'$ at
$s'$, $I=I'$, and $S=S'$. We have $P|_{<L} = [ T' S] =[T'S']=P'$,
proving \eqref{eq:Presprime}. $Q'$ records the insertion of
$I'w_0=Iw_0$ into $T'$. $Q$ records the insertion of the same word
into $T'$ except that there is a large number $L$ present at the
cell $s'$. Then $Q$ is given by the internal insertion $I_{s'}(Q')$.
\item Suppose $x=\bar{r}$. Then $T$ is obtained by the reverse row
insertion on $T'$ at the cell $s'$, ejecting the letter $i$, say. In
particular $[T']=[T i]$. Also $I=I'(i,L)$, and $S$ is obtained from
$S'$ by the column insertion of $i$ followed by the adjoining of $L$
at the end of the column just right of the end of the bumping path
of the column insertion. It follows that $S|_{<L}=[i S']$. We have
$P|_{<L}=[T (i S')]=[(Ti)S']=[T'S']=P'$, proving
\eqref{eq:Presprime}. The biword $Iw_0$ is obtained from $I'w_0$ by
adding biletters $(i,L^*)$ and $(L,i^*)$. The biletter $(i,L^*)$ is
the first in $Iw_0$ to be inserted into $T$. This produces the
tableau $T'=[Ti]$, so that $L^*$ goes in the cell $s'$ of $Q$.
Thereafter it is the same as inserting $I' w_0$ into $T'$, except
that in the middle there is a biletter $(L,i^*)$ containing the
largest value $L$. It follows that $Q$ is obtained from $Q'$ by the
external insertion of $i^*$ and putting the smallest letter $L^*$ at
the cell $s'$ (which is on the inside of the skew shape of $Q'$).
\item Suppose $x=\nn$. Then $T=T'$, $I$ is obtained from $I'$ by
adding the fixed point $L$, and  $S$ is obtained from $S'$ by
adjoining $L$ at the end of the first column, that is, $S=[L S']$.
Therefore $P|_{<L}=[T S']=[T'S']=P'$, proving \eqref{eq:Presprime}.
$Iw_0$ is obtained from $I'w_0$ by putting the biletter $(L,L^*)$ at
the beginning. Since $T=T'$ the only difference between the
computations of $Q$ and $Q'$ is that one first row inserts $L$ and
records by $L^*$. Therefore $Q$ is obtained from $Q'$ by the
external insertion of $L^*$.
\end{enumerate}

\subsection{$DDF=VXR$ and the proof of $X=K$}
\begin{prop} \label{pp:PQV=D} $(P,Q,T)=(\Pt,\Qt,\Tt)$.
\end{prop}
\begin{proof}
By induction $(P',Q',T')=(\Pt',\Qt',\Tt')$ where the primed elements
correspond to the shorter path $b'$ where $b=x\otimes b'$. By the
analysis in sections \ref{ss:incDDF} and \ref{ss:incVXR}, the
changes $(P',Q',T')\rightarrow (P,Q,T)$ and
$(\Pt',\Qt',\Tt')\rightarrow (\Pt,\Qt,\Tt)$ are the same.
\end{proof}

We now finish the proof of Theorem \ref{th:KX}. It suffices to prove
Theorem \ref{th:VXR}. By Proposition \ref{pp:nu} it suffices to
prove Theorem \ref{th:VXR} for the case $\nu=1^L$. The VXR map
$\nu_{1^L}$ is bijective since it agrees with the DDF map, which is
known to be a bijection (Propositions \ref{pp:PQV=D} and
\ref{pp:DDFbij}). Finally, \eqref{eq:bijstat} holds by Proposition
\ref{pp:stat}. This concludes the proof of Theorem \ref{th:KX}.

\section{Schensted and DDF according to Fomin and Roby}
\label{sec:Roby} Although not needed for the proofs, we include an
interpretation of the DDF map using Fomin's growth diagrams,
extending Roby's interpretation of the original DDF bijection
\cite{R95}. This is insightful, as one sees the local coenergy
calculations mirrored in this computation of the DDF map.

\subsection{Fomin's growth diagrams}
Fomin observed that various Schensted algorithms may be understood
solely in terms of chains of partitions \cite{F94,F95}. Indeed, many
of the properties of these algorithms are most clearly seen from
this perspective.

An injective tableau $T$ of shape $\la/\mu$ in the alphabet
$\{1,2,\dotsc,N\}$ may be interpreted as the sequence of partitions
$\mu=\la^{(0)}\subset\la^{(1)}\subset\dotsm\subset \la^{(N)}=\la$
where each $\la^{(i)}=\la^{(i-1)}$ if $i$ does not appear in $T$ and
otherwise $\la^{(i)}$ is obtained from $\la^{(i-1)}$ by adding the
cell of $T$ containing $i$. Or, the outer shape of each $T|_{\le i}$
is $\la^{(i)}$, where $T|_{\le i}$ is the skew tableau given by the
restriction of $T$ to the letters of value at most $i$. Note that
the letter $i$ is associated with the transition from the shape
$\la^{(i-1)}$ to $\la^{(i)}$.

Consider the Sagan-Stanley skew RS map $(T,U,w)\rightarrow (P,Q)$
given in subsection \ref{ss:skew}. There is an associated
\textbf{shape array} $(\la_a^b)$, a two-dimensional array of
partitions with row index set of size $|B|+1$ and column index set
of size $|A|+1$. Let us add a smallest letter $a_0$ to $A$ and $b_0$
to $B$. The array elements $\la_a^b$ will be indexed using these
enlarged alphabets. The letters of $A$ (resp. $B$) naturally label
the transitions from one row (resp. column) to the next. In other
words, the pair $(a,b)\in A\times B$ naturally labels the point in
the middle of a two-by-two submatrix, whose lower right entry is
$\la_a^b$. In particular a pair $(w(b),b)$ for $b\in B'$, is
naturally associated with such a two-by-two submatrix.

Let $w|_{\le a}^{\le b}$ be the restriction of $w$ to the positions
at most $a$ and values at most $b$. In other words, viewing $w$ as a
set of ordered pairs $(w(b),b)$, select only those pairs $(a',b')$
such that $a'\le a$ and $b'\le b$.

The entry $\la_a^b$ is defined to be the shape of the bumping
tableau associated with the data $(T|_{\le a},U|_{\le b},w|_{\le
a}^{\le b})$ under skew RS. Then the left and top borders of the
shape array give $T$ and $U$ respectively, and the right and bottom
borders give $P$ and $Q$ respectively.

Fomin's observation is that given $(T,U,w)$, one may recover the
entire shape array using a local rule. Set the left and top borders
of the shape array to be $T$ and $U$ respectively, with all other
entries to be determined.

Consider a two-by-two subarray of partitions
\begin{equation} \label{eq:2by2}
\young(\mu\nu,\rho\sigma)
\end{equation}
indexed by $(a,b)$, in which all but $\sigma$ have already been
determined. Let $s$ and $t$ be the cells (possibly empty) that must
be adjoined to $\mu$ to obtain $\nu$ and $\rho$ respectively. Then
$\sigma$ is defined by the following rules.
\begin{enumerate}
\item If $s=t$ is not empty then $\sigma$ is obtained from
$\nu=\rho$ by adjoining a cell to the $(r+1)$-st row where $s=t$ is
in the $r$-th row of $\nu=\rho$.
\item If $s=t$ is empty and $w(b)=a$ then $\sigma$ is
obtained from $\mu=\nu=\rho$ by adjoining a cell to the first row.
\item Otherwise $\sigma=\nu\cup\rho$.
\end{enumerate}
Applying this local rule one may uniquely recover the entire array
and read off $P$ and $Q$ at the end.

Or, starting with $P$ and $Q$ at the right and bottom borders with
the rest of the shape array to be determined, one may apply the
reverse of the above local rule. Suppose in \eqref{eq:2by2} (indexed
by the pair $(a,b)$) one knows all of the partitions but $\mu$. Then
$\mu$ can be recovered, as well as whether or not $w(b)=a$.

Let $s'$ and $t'$ be the cells (possibly empty) which must be
removed from $\sigma$ to obtain $\rho$ and $\nu$ respectively.
\begin{enumerate}
\item If $s'=t'$ is nonempty and $s'=t'$ is in the first row of
$\sigma$, then add the pair $(a,b)$ to $w$ (that is, extend the
bijection $w$ by declaring that $w(b)=a$) and set $\mu$ equal to
$\nu=\rho$. We call such a two-by-two array \textbf{special}.
\item If $s'=t'$ is nonempty and $s'=t'$ is in the $(r+1)$-st row of
$\sigma$ with $r\ge1$, then $\mu$ is obtained from $\nu=\rho$ by
removing a cell in the $r$-th row.
\item Otherwise $\mu=\nu\cap\rho$.
\end{enumerate}
Using the reverse local rule the rest of the shape array may be
determined, yielding the right and top borders $T$ and $U$ and the
bijection $w$.

\begin{theorem} \cite{SS90} \label{th:RSinv} The skew RS map
defined in subsection \ref{ss:skew} has the involution property,
that is, $(T,U,w)\mapsto (P,Q)$ if and only if $(U,T,w^{-1})\mapsto
(Q,P)$.
\end{theorem}
\begin{proof} Fomin's local rule is invariant under the transposition of
the two-by-two matrix \eqref{eq:2by2}.
\end{proof}

\subsection{Roby on DDF}
Roby \cite{R95} observed that a skew extension \cite{DS95} of the
DDF bijection could be computed using Fomin's local rule. We only
require the nonskew case.

Let $\lad$ be an oscillating tableau of shape $\la$ and length $L$
and $(T,I)$ its image under the DDF bijection. The pair $(T,I)$ may
be recovered as follows. In the skew RS algorithm, set both
alphabets $A$ and $B$ to be $\{1,2,\dotsc,L\}$. Start with an
unfilled $(L+1)\times (L+1)$ shape array. Put the oscillating
tableau on the \textit{antidiagonal}, starting from the upper right
and proceeding to the lower left. Use the reverse of Fomin's local
rule to fill in the half of the diagram above the antidiagonal. To
start the process one must be able to compute $\mu$ in the reverse
local rule for two-by-two subarrays \eqref{eq:2by2} such that $\nu$
and $\rho$ are consecutive shapes in the oscillating tableau,
without knowledge of $\sigma$. But consecutive shapes in oscillating
tableaux are never equal, so one does not need $\sigma$ to determine
$\mu$. Filling in the triangle of the shape array above the
antidiagonal, the left border recovers the injective tableau $T$.
The top border consists of empty partitions since the oscillating
tableau starts with the empty partition. One also obtains a set of
special positions $(a,i)$ above the antidiagonal, which yield the
permutation data: $I$ is the product of the transpositions of the
form $(a,\hat{i})=(a,L+1-i)$.

This is readily extended for Motzkin tableaux. In this case there is
ambiguity on how to compute the reverse local rule when consecutive
shapes in the Motzkin tableau are equal. We resolve this ambiguity
by declaring that if $\la^{(i-1)}=\la^{(i)}$ then the corresponding
position $(i,\hat{i})$ is special. In other words, for a two-by-two
subarray \eqref{eq:2by2} on the antidiagonal, if $\nu=\rho$ we
declare that $\sigma$ is obtained from $\nu=\rho$ by adding a cell
to the first row, and set $\mu$ equal to $\nu=\rho$. Then there is a
unique way to fill in the triangular part of the shape array above
the antidiagonal. Again $T$ is the left border and the top border
consists of empty partitions. As before, the special positions
strictly above the antidiagonal give the transpositions in the
desired involution $I$, while the fixed points $i$ of $I$ are given
by the special positions $(i,\hat{i})$ on the antidiagonal.

See Figure \ref{fig:Mshape} for the shape array associated with the
Motzkin tableau given in Figure \ref{fig:mDDF}, in which the special
positions are indicated by the symbol $\otimes$.

\begin{figure}
\Yboxdim{3pt}\Yinterspace{1pt}
\begin{equation*}
\begin{array}{|c||ccccccccccccccc|} \hline
\raisebox{5pt}{\vphantom{1}} && \hg && \hf && \he && \hd && \hc && \hb && \ha &\\
 \hline \hline
 & \dt && \dt && \dt && \dt && \dt && \dt && \dt && \dt \\ 1 &&\otimes&&&&&&&&&&&&& \\
 & \dt &&  \yng(1) &&\yng(1) &&\yng(1) &&\yng(1) &&\yng(1) &&\yng(1)
&&\yng(1) \\ 2 &&&&&&&&\otimes&&&&&&& \\
 & \dt && \yng(1) &&\yng(1) &&\yng(1) && \yng(2) &&\yng(2)
&&\yng(2) &&\yng(2) \\ 3 &&&&&&&&&&&&&&&\\
 & \yng(1) && \yng(1,1) &&\yng(1,1) &&\yng(1,1) &&\yng(2,1) &&\yng(2,1)
&&\yng(2,1) &&\yng(2,1) \\ 4 &&&&&&&&&&&&\otimes &&&\\
 & \yng(1) && \yng(1,1) &&\yng(1,1) &&\yng(1,1) &&\yng(2,1) &&\yng(2,1)
&&\yng(3,1) &&\yng(3,1) \\ 5 &&&&&&\otimes &&&&&&&&&\\
 & \yng(1) && \yng(1,1) &&\yng(1,1) &&\yng(2,1) &&\yng(2,2) &&\yng(2,2) &&\yng(3,2) &&\yng(3,2)
 \\ 6 &&&&&&&&&&&&&&&\\
 & \yng(2) && \yng(2,1) &&\yng(2,1) &&\yng(2,2) &&\yng(2,2,1) &&\yng(2,2,1)
&&\yng(3,2,1) &&\yng(3,2,1) \\ 7 &&&&&&&&&&&&&&\otimes&\\
 & \yng(2) && \yng(2,1) &&\yng(2,1) &&\yng(2,2) &&\yng(2,2,1) &&\yng(2,2,1)
&&\yng(3,2,1) &&\yng(4,2,1) \\ \hline
\end{array}
\end{equation*}
\Yboxdim{13pt}\Yinterspace{1pt} \caption{Shape array for a Motzkin
path} \label{fig:Mshape}
\end{figure}

\subsection{The lower triangle of the shape array}
For our application the lower triangle of the shape array is
crucial. Starting with the Motzkin tableau of shape $\la$, the upper
triangle and special positions (and therefore $T$ and $I$) are
determined as before. To fill out the positions in the lower
triangle below the antidiagonal it is necessary to know which of
those positions are special. We define the special positions below
the antidiagonal to be those whose mirror images across the
antidiagonal, are special. This is equivalent to defining the total
set of special positions to be given by the bijection $Iw_0$ defined
in section \ref{ss:DDFins}.

The triangle below the antidiagonal may be filled using Fomin's
local rule. At the end, let $P$ and $Q$ be the right and bottom
borders. $P$ is a standard tableau of some shape, say $\tau$, and
$Q$ is an injective skew tableau of shape $\tau/\la$. Then the skew
RS algorithm gives
\begin{equation} \label{eq:DDFSS}
(T,\nn,Iw_0)\mapsto (P,Q).
\end{equation}
This coincides with the $P$ tableau defined in \eqref{eq:DDFP}.

\begin{example} \label{ex:RobyDDF} Figure \ref{fig:mDDF} contains
the data $I w_0,P,Q$ for each of the Motzkin tableaux associated
with right factors $b_i\dotsm b_2 b_1$ of $b$ in Example
\ref{ex:hwv}.
\end{example}

\subsection{Ghost of coenergy in DDF}
Let $b\in P_\se(1^L,\la)$ and $\dck\otimes c$ be as in
\eqref{eq:dcdef}. Remarkably, the computation of the coenergy
$\Db(\dck\otimes c)$ in section \ref{ss:grading} can be seen clearly
in the Fomin shape array of the DDF map applied to $b$. Compare the
underlined positions in the computation of $R_+$ given in Example
\ref{ex:dckc}, with the special positions in the shape array for $b$
in Figure \ref{fig:Mshape} that are on or southeast of the main
antidiagonal. They coincide if the shape array is rotated 45 degrees
clockwise.

\section{The missing cases}
\label{sec:end} Our methods don't apply when $\se=\vdom\,$. Although
the paths for types $\se=\vdom$ and $\se=\hdom$ for tensor powers of
$B_\se$ are in an easy bijection given by ``transposing shapes", the
coenergy functions are quite different in nature. This case requires
a new idea.

We also expect that the VXR map will succeed for any tensor product
of KR crystals for types $\se\in\onetwoset$. To do this one must
construct the KR modules $W^{(r)}_s$ and crystals $B^{r,s}$ in
general and show, as conjectured in \cite{OSS03} that the virtual
crystal methods work for these crystals. For the case of tensor
products of KR crystals of the form $B^{r,1}$ the virtual crystals
have already been constructed \cite{OSS03}. Our methods should work
in this case, using intrinsic energy rather than coenergy and
compatible splitting maps that remove a box from a column. Note that
this case is not the transpose of the case proven in Theorem
\ref{th:KX}; the asymmetry is seen in the classical decomposition of
the KR modules (the $U'_q(C_n^{(1)})$-crystal $B^{r,1}$ is
isomorphic to $B(\omega_r)$ as $U_q(C_n)$-crystals) and in the $K$
formula in the case $\hdom\,$ where the horizontal dominoes do not
get transposed.


\begin{thebibliography}{xxx}

\bibitem{AK97} T. Akasaka and M. Kashiwara,
Finite-dimensional representations of quantum affine algebras, Publ.
Res. Inst. Math. Sci. \textbf{33} (1997) 839--867.

\bibitem{B74} W. Burge, Four correspondences between graphs and
generalized Young tableaux, J. Combinatorial Theory Ser. A
\textbf{17} (1974), 12--30.

\bibitem{Ch01} V. Chari, On the fermionic formula and the
Kirillov-Reshetikhin conjecture, Internat. Math. Res. Notices
2001, 629--654.

\bibitem{CP95} V. Chari and A. Pressley, Minimal affinizations of
representations of quantum groups: the nonsimply-laced case. Lett.
Math. Phys. \textbf{35} (1995) 99--114.

\bibitem{CK00} V. Chari and M. Kleber, Symmetric functions and
representations of quantum affine algebras, Recent developments in
infinite-dimensional Lie algebras and conformal field theory
(Charlottesville, VA, 2000), 27--45, Contemp. Math., \textbf{297},
Amer. Math. Soc., Providence, RI, 2002.

\bibitem{DDF88} M.-P.~Delest, S.~Dulucq and L.~Favreau,
An analogue to Robinson-Schensted correspondence for oscillating
tableaux, S\'eminaire Lotharingien de Combinatoire, B20b (1988).

\bibitem{DS95} S. Dulucq and B. Sagan, La correspondance
de Robinson-Schensted pour les tableaux oscillants gauches, Formal
power series and algebraic combinatorics (Montreal, PQ, 1992).
Discrete Math. \textbf{139} (1995), no. 1-3, 129--142.

\bibitem{F94} S. Fomin, Duality of graded graphs, J. Algebraic Combin.
\textbf{3} (1994),  no. 4, 357--404.

\bibitem{F95} S. Fomin, Schensted algorithms for dual graded graphs, J.
Algebraic Combin. \textbf{4} (1995),  no. 1, 5--45.

\bibitem{FG93} S. Fomin and C. Greene, A Littlewood-Richardson miscellany,
European J. Combin. \textbf{14} (1993), no. 3, 191--212.

\bibitem{Ful97} W. Fulton, Young tableaux, with applications to representation
theory and geometry. London Mathematical Society Student Texts,
\textbf{35}, Cambridge University Press, Cambridge, 1997.

\bibitem{Her05} D. Hernandez, The Kirillov-Reshetikhin conjecture and solutions of
T-systems, arXiv:math.QA/0501202.

\bibitem{HKOT00} G.~Hatayama, A.~Kuniba, M.~Okado,
T.~Takagi, Combinatorial $R$ matrices for a family of crystals:
$C\sb n\sp {(1)}$ and $A\sb {2n-1}\sp {(2)}$ cases, Physical
combinatorics (Kyoto, 1999), 105--139, Progr. Math., 191, Birkhäuser
Boston, Boston, MA, 2000.

\bibitem{HKOT02} G.~Hatayama, A.~Kuniba, M.~Okado,
T.~Takagi, Combinatorial $R$ matrices for a family of crystals:
$B\sp {(1)}\sb n,D\sp {(1)}\sb n,A\sp {(2)}\sb {2n}$, and $D\sp
{(2)}\sb {n+1}$ cases, J. Algebra 247 (2002) 577--615.

\bibitem{HKOTT01} G.~Hatayama, A.~Kuniba, M.~Okado, T.~Takagi, and Z.~Tsuboi,
Paths, Crystals and Fermionic Formulae, arXiv:math.QA/0102113.

\bibitem{HKOTY99} G.~Hatayama, A.~Kuniba, M.~Okado, T.~Takagi, and
Y.~Yamada, Remarks on fermionic formula, in Recent developments in
quantum affine algebras and related topics (Raleigh, NC, 1998),
243--291, Contemp. Math., \textbf{248}, Amer. Math. Soc.,
Providence, RI, 1999.

\bibitem{JMO00} N.~Jing, K.~C.~Misra, and M.~Okado,
$q$-wedge modules for quantized enveloping algebras of classical
type, J. Algebra \textbf{230} (2000), no. 2, 518--539.

\bibitem{Kac90} V. Kac, Infinite-dimensional Lie algebras, third
edition, Cambridge University Press, Cambridge, 1990.

\bibitem{Kas95} M. Kashiwara, On crystal bases. Representations of groups (Banff, AB, 1994),
155--197, CMS Conf. Proc., \textbf{16}, Amer. Math. Soc.,
Providence, RI, 1995.

\bibitem{Kas03} M. Kashiwara,
Level zero fundamental representations over quantized affine
algebras and Demazure modules, arXiv:math.QA/0309142.

\bibitem{KKM94} S.~J.~Kang, M.~Kashiwara, and K.~C.~Misra,
Crystal bases of Verma modules for quantum affine Lie algebras.
Compositio Math. \textbf{92} (1994),  no. 3, 299--325.

\bibitem{KKR88} S.~.V.~Kerov, A.~N.~Kirillov, and
N.~Yu.~Reshetikhin, Combinatorics, the Bethe ansatz and
representations of the symmetric group, J. Soviet Math.
\textbf{41} (1988), no. 2, 916--924.

\bibitem{Kl02} M. Kleber, Embeddings of Schur functions into types
$B/C/D$. J. Algebra \textbf{247} (2002), no. 2, 452--466.

\bibitem{KMN92} S.-J.~Kang, M.~Kashiwara, K.~C.~Misra, T.~Miwa, T.~Nakashima, and
A.~Nakayashiki, Affine crystals and vertex models, Infinite
analysis, Part A, B (Kyoto, 1991), 449--484, Adv. Ser. Math.
Phys., \textbf{16}, World Sci. Publishing, River Edge, NJ, 1992.

\bibitem{KMN92a} S.-J.~Kang, M.~Kashiwara, K.~C.~Misra, T.~Miwa, T.~Nakashima, and
A.~Nakayashiki, Perfect crystals of quantum affine Lie algebras.
Duke Math. J.  \textbf{68} (1992),  no. 3, 499--607.

\bibitem{KN94} M. Kashiwara and T. Nakashima,
Crystal graphs for representations of the $q$-analogue of classical
Lie algebras, J. Algebra \textbf{165} (1994) 295--345.

\bibitem{Kn70} D. E. Knuth, Permutations, matrices, and generalized Young
tableaux, Pacific J. Math. \textbf{34} (1970) 709--727.

\bibitem{Ko99} Y. Koga, Level one perfect crystals for $B_n^{(1)}$, $C_n^{(1)}$,
and $D_n^{(1)}$, J. Algebra \textbf{217} (1999), no. 1, 312--334.

\bibitem{KR88} A.~N.~Kirillov, and
N.~Yu.~Reshetikhin, The Bethe ansatz and the combinatorics of
Young tableaux, J. Soviet Math. \textbf{41} (1988), no. 2,
925--955.

\bibitem{KR90} A. N. Kirillov and N.~Yu.~Reshetikhin, Representations of Yangians
and multiplicities of the inclusion of the irreducible components
of the tensor product of representations of simple Lie algebras,
J. Soviet Math. \textbf{52} (1990),  no. 3, 3156--3164.

\bibitem{KSS02} A.~N.~Kirillov, A. Schilling, and M. Shimozono,
A bijection between Littlewood-Richardson tableaux and rigged
configurations, Selecta Math. (N.S.) \textbf{8} (2002),  no. 1,
67--135.

\bibitem{L91} A. Lascoux, Cyclic permutations on words, tableaux
and harmonic polynomials.  Proceedings of the Hyderabad Conference
on Algebraic Groups (Hyderabad, 1989), 323--347, Manoj Prakashan,
Madras, 1991.

\bibitem{Le04} C. Lecouvey, A duality between $q$-multiplicities in tensor
products and $q$-multiplicities of weights for the root systems $B$,
$C$ or $D$, preprint, arXiv:math.RT/0407522.

\bibitem{Le04a} C. Lecouvey, Branching rules, Kostka-Foulkes polynomials and
$q$-multiplicities in tensor product for the root systems $B_n$,
$C_n$ and $D_{n}$, preprint arXiv:math.CO/0412548.

\bibitem{LR34} D. E. Littlewood and A. R. Richardson, Group
characters and algebra, Philos. Trans. R. Soc., A, \textbf{233}
(1934) 99-141.

\bibitem{LS81} A. Lascoux and M. P. Sch\"utzenberger, Le Mono\"ide
Plaxique, Noncommutative structures in algebra and geometric
combinatorics (Naples, 1978), pp. 129--156, Quad. Ricerca Sci.,
\textbf{109}, CNR, Rome, 1981.

\bibitem{Lu83} G. Lusztig, Singularities, character formulas, and a
$q$-analog of weight multiplicities, Analysis and topology on
singular spaces, II, III (Luminy, 1981), 208--229, Ast\'erisque,
101-102, Soc. Math. France, Paris, 1983.

\bibitem{Mac95} I. G. Macdonald,
Symmetric functions and Hall polynomials, Second edition, With
contributions by A. Zelevinsky, Oxford Mathematical Monographs,
Oxford Science Publications, The Clarendon Press, Oxford
University Press, New York, 1995.

\bibitem{NY97} A. Nakayashiki and Y. Yamada, Kostka polynomials and
energy functions in solvable lattice models, Selecta Math. (N.S.)
\textbf{3} (1997), no. 4, 547--599.

\bibitem{OSS03} M. Okado, A. Schilling, and M. Shimozono, Virtual
crystals and fermionic formulas of type $D^{(2)}_{n+1}$,
$A^{(2)}_{2n}$, and $C^{(1)}_n$, Represent. Theory  \textbf{7}  (2003),
101--163 (electronic).

\bibitem{OSS03a} M. Okado, A. Schilling, and M. Shimozono,
A crystal to rigged configuration bijection for nonexceptional
affine algebras, ``Algebraic Combinatorics and Quantum Groups",
Edited by N. Jing, World Scientific (2003), 85--124.

\bibitem{OSS03b} M. Okado, A. Schilling, and M. Shimozono,
Virtual crystals and Kleber's algorithm, Commun. Math. Phys.
\textbf{238} (2003) 187-209.

\bibitem{R95} T. Roby, The connection between the Robinson-Schensted
correspondence for skew oscillating tableaux and graded graphs,
Formal power series and algebraic combinatorics (Montreal, PQ,
1992), Discrete Math. \textbf{139}  (1995),  no. 1-3, 481--485.

\bibitem{Sc61} C. Schensted, Longest increasing and decreasing
subsequences, Canad. J. Math. \textbf{13} (1961) 179--191.

\bibitem{Sci04} A. Schilling,
A bijection between type $D_n^{(1)}$ crystals and rigged
configurations, preprint arXiv:math.QA/0406248.

\bibitem{ScS04} A. Schilling and M. Shimozono, $X=M$ for
symmetric powers, preprint arXiv:math.QA/0412376.

\bibitem{ScWa99} A. Schilling and S. Ole Warnaar,
Inhomogeneous lattice paths, generalized Kostka polynomials and
$A_{n-1}$ supernomials, Comm. Math. Phys. \textbf{202} (1999)
359--401.

\bibitem{S02} M. Shimozono, Affine type A crystal structure
on tensor products of rectangles, Demazure characters, and nilpotent
varieties.  J. Algebraic Combin. \textbf{15} (2002) 151--187.

\bibitem{S02a} M. Shimozono, Multi-atoms and monotonicity of
generalized Kostka polynomials, European J. Combin. \textbf{22}
(2001) 395--414.

\bibitem{SS90} B. Sagan and R. P. Stanley, Robinson-Schensted algorithms
for skew tableaux, J. Combin. Theory Ser. A \textbf{55}  (1990),
no. 2, 161--193.

\bibitem{St01} J. Stembridge, Multiplicity-free products of Schur
functions, Ann. Comb. \textbf{5} (2001), 113--121.

\bibitem{SZ04} M. Shimozono and M. Zabrocki, Deformed universal characters
for classical and affine algebras, preprint arXiv:math.CO/0404288.

\bibitem{Y98} S. Yamane, Perfect crystals of $U_q(G_2^{(1)})$,
J. Algebra \textbf{210} (1998),  no. 2, 440--486.

\end{thebibliography}
\end{document}